\documentclass[english]{IEEEtran}
\usepackage[T1]{fontenc}
\usepackage{babel}
\usepackage{float}
\usepackage{amsthm}
\usepackage{amsmath}
\usepackage{amssymb}
\usepackage{graphicx}
\usepackage[unicode=true,
 bookmarks=true,bookmarksnumbered=true,bookmarksopen=true,bookmarksopenlevel=1,
 breaklinks=false,pdfborder={0 0 0},backref=false,colorlinks=false]
 {hyperref}
\hypersetup{pdftitle={Your Title},
 pdfauthor={Your Name},
 pdfpagelayout=OneColumn, pdfnewwindow=true, pdfstartview=XYZ, plainpages=false}
\usepackage{breakurl}

\makeatletter

\providecommand{\tabularnewline}{\\}
\floatstyle{ruled}
\newfloat{algorithm}{tbp}{loa}
\providecommand{\algorithmname}{Algorithm}
\floatname{algorithm}{\protect\algorithmname}

\theoremstyle{plain}
\newtheorem{thm}{\protect\theoremname}
\theoremstyle{plain}
\newtheorem{lem}[thm]{\protect\lemmaname}
 \usepackage{algolyx}

\ifCLASSOPTIONcompsoc
\usepackage[caption=false,font=normalsize,labelfont=sf,textfont=sf]{subfig}
\else
\usepackage[caption=false,font=footnotesize]{subfig}
\fi

\usepackage{cite}

\makeatletter
\def\thmhead@plain#1#2#3{%
  \thmname{#1}\thmnumber{\@ifnotempty{#1}{ }\@upn{#2}}%
  \thmnote{ {\the\thm@notefont#3}}}
\let\thmhead\thmhead@plain
\makeatother

\makeatother

\providecommand{\lemmaname}{Lemma}
\providecommand{\theoremname}{Theorem}

\begin{document}

\title{Sequence Set Design With Good Correlation Properties via Majorization-Minimization}

\author{Junxiao Song, Prabhu Babu, and Daniel P. Palomar, \IEEEmembership{Fellow, IEEE}\thanks{This work was supported by the Hong Kong RGC 16206315 research grant.
Junxiao Song, Prabhu Babu, and Daniel P. Palomar are with the Hong
Kong University of Science and Technology (HKUST), Hong Kong. E-mail:
\{jsong, eeprabhubabu, palomar\}@ust.hk.}}
\maketitle
\begin{abstract}
Sets of sequences with good correlation properties are desired in
many active sensing and communication systems, e.g., multiple-input\textendash multiple-output
(MIMO) radar systems and code-division multiple-access (CDMA) cellular
systems. In this paper, we consider the problems of designing complementary
sets of sequences (CSS) and also sequence sets with both good auto-
and cross-correlation properties. Algorithms based on the general
majorization-minimization method are developed to tackle the optimization
problems arising from the sequence set design problems. All the proposed
algorithms can be implemented by means of the fast Fourier transform
(FFT) and thus are computationally efficient and capable of designing
sets of very long sequences. A number of numerical examples are provided
to demonstrate the performance of the proposed algorithms.\end{abstract}


\section{Introduction}

Sequences with good correlation properties play an important role
in many active sensing and communication systems\cite{he2012waveform,levanon2004radar}.
The design of a single sequence with good autocorrelation properties
(e.g., small autocorrelation sidelobes) has been studied extensively,
e.g., see \cite{stoica2009new,MISL,WISL_song2015} and the references
therein. In this paper, we focus on the design of sets of sequences
with good correlation properties. We consider both the design of complementary
sets of sequences (CSS) and the design of sequence sets with good
auto- and cross-correlation properties. In addition, in order to avoid
non-linear side effects and make full use of the transmission power
available in the system, we restrict our design to unimodular sequences.

Let $\{\mathbf{x}_{m}\}_{m=1}^{M}$ denote a set of $M$ complex unimodular
sequences each of length $N$, i.e., $\mathbf{x}_{m}=[x_{m}(1),\ldots,x_{m}(N)]^{T}$,
$m=1,\ldots,M$. Then the aperiodic cross-correlation of $\mathbf{x}_{i}$
and $\mathbf{x}_{j}$ at lag $k$ is defined as 
\begin{eqnarray}
r_{i,j}(k) & = & \sum_{n=1}^{N-k}x_{i}(n+k)x_{j}^{*}(n)=r_{j,i}^{*}(-k),\nonumber \\
 &  & i,j=1,\ldots,M,k=1-N,\ldots,N-1.\label{eq:cross_corr-1}
\end{eqnarray}
When $i=j,$ \eqref{eq:cross_corr-1} reduces to the autocorrelation
of $\mathbf{x}_{i}$.

The motivation of CSS design comes from the difficulties in designing
a single unimodular sequence with impulse-like autocorrelation. For
instance, it can be easily observed that the autocorrelation sidelobe
at lag $N-1$ of a unimodular sequence is always equal to 1, no matter
how we design the sequence. The difficulties have encouraged researchers
to consider the idea of CSS, and the set of sequences $\{\mathbf{x}_{m}\}_{m=1}^{M}$
is called complementary if and only if the autocorrelations of $\{\mathbf{x}_{m}\}$
sum up to zero at any out-of-phase lag, i.e., 
\begin{equation}
\sum_{m=1}^{M}r_{m,m}(k)=0,\,1\leq\left|k\right|\leq N-1.
\end{equation}
CSS have been applied in many active sensing and communication systems,
for instance, multiple-input\textendash multiple-output (MIMO) radars
\cite{compleSet_Searle2008}, radar pulse compression \cite{compleSet_Levanon2009},
orthogonal frequency-division multiplexing (OFDM) \cite{compleSet_Schmidt2007},
ultra wide-band (UWB) communications \cite{compleSet_Garcia2010},
code-division multiple-access (CDMA) \cite{compleSet_Tseng2000},
and channel estimation \cite{compleSet_Spasojevic2001}. Owing to
the practical importance, a lot of effort has been devoted to the
construction of CSS. The majority of research results on CSS at the
early stage have been concerned with the analytical construction of
CSS for restricted sequence length $N$ and set cardinality $M$.
More recently, computational methods have also been proposed for the
design of CSS, see \cite{SequenceSet_Mojtaba2013} for example. In
contrast to analytical constructions, computational methods are more
flexible in the sense that they do not impose any restriction on the
length of sequences or the set cardinality. 

In CSS design, only the autocorrelation properties of the sequences
have been considered. But some applications require a set of sequences
with not only good autocorrelation properties but also good cross-correlations
among the sequences, for example, in CDMA cellular networks or in
MIMO radar systems. Good autocorrelation indicates that a sequence
is nearly uncorrelated with its own time-shifted versions, while good
cross-correlation means that any sequence is nearly uncorrelated with
all other time-shifted sequences. Good correlation properties in the
above sense ensure that matched filters at the receiver end can easily
separate the users in a CDMA system \cite{lowCorrSet_Oppermann1997}
or extract the signals backscattered from the range of interest while
attenuating signals backscattered from other ranges in MIMO radar
\cite{SequenceSet_He2009}. 

Extending the approaches in \cite{WISL_song2015}, we present in this
paper several new algorithms for the design of complementary sets
of sequences and sequence sets with both good auto- and cross-correlation
properties. The sequence set design problems are first formulated
as optimization problems and they include the single sequence design
problems considered in \cite{MISL,WISL_song2015} as special cases.
Then several efficient algorithms are developed based on the general
majorization-minimization (MM) method via successively majorizing
the objective functions twice. All the proposed algorithms can be
implemented by means of the fast Fourier transform (FFT) and are thus
very efficient in practice. The convergence properties and an acceleration
scheme, which can be used to further accelerate the proposed MM algorithms,
are also briefly discussed.

The remaining sections of the paper are organized as follows. In Section
\ref{sec:Problem-Formulation}, the problem formulations are presented.
In Section \ref{sec:Seq-Design-viaMM}, an MM algorithm is derived
for the CSS design problem, followed by the derivations of two MM
algorithms for designing sequence sets with good auto- and cross-correlations
in Sections \ref{sec:Set_Auto_Cross_Wei} and \ref{sec:Adaptive-MM},
respectively. Convergence analysis and an acceleration scheme are
introduced in Section \ref{sec:Convergence-Acc}. Finally, Section
\ref{sec:Numerical-Experiments} presents some numerical results,
and the conclusions are given in Section \ref{sec:Conclusion}.

\emph{Notation}: Boldface upper case letters denote matrices, boldface
lower case letters denote column vectors, and italics denote scalars.
$\mathbb{R}$ and $\mathbb{C}$ denote the real field and the complex
field, respectively. $\mathrm{Re}(\cdot)$ and $\mathrm{Im}(\cdot)$
denote the real and imaginary part, respectively. ${\rm arg}(\cdot)$
denotes the phase of a complex number. The superscripts $(\cdot)^{T}$,
$(\cdot)^{*}$ and $(\cdot)^{H}$ denote transpose, complex conjugate,
and conjugate transpose, respectively. $X_{i,j}$ denotes the (\emph{i}-th,
\emph{j}-th) element of matrix $\mathbf{X}$ and $x_{i}$ ($x(i)$)
denotes the \emph{i}-th element of vector $\mathbf{x}$. $\mathbf{X}_{i,:}$
denotes the \emph{i}-th row of matrix $\mathbf{X}$, $\mathbf{X}_{:,j}$
denotes the \emph{j}-th column of matrix $\mathbf{X}$, and $\mathbf{X}_{i:j,k:l}$
denotes the submatrix of $\mathbf{X}$ from $X_{i,k}$ to $X_{j,l}$.
$\circ$ denotes the Hadamard product. $\otimes$ denotes the Kronecker
product. $\mathrm{Tr}(\cdot)$ denotes the trace of a matrix. ${\rm diag}(\mathbf{X})$
is a column vector consisting of all the diagonal elements of $\mathbf{X}$.
${\rm Diag}(\mathbf{x})$ is a diagonal matrix formed with $\mathbf{x}$
as its principal diagonal. ${\rm vec}(\mathbf{X})$ is a column vector
consisting of all the columns of $\mathbf{X}$ stacked. $\mathbf{I}_{n}$
denotes an $n\times n$ identity matrix.

\section{Problem Formulation and MM Primer\label{sec:Problem-Formulation}}

The problems of interest in this paper are the design of complementary
sets of sequences (CSS) and the design of sequence sets with good
auto- and cross-correlation properties. In the following, we first
provide criteria to measure the complementarity of a sequence set
and also the goodness of auto- and cross-correlation properties respectively,
and then formulate the sequence set design problems as optimization
problems. The MM method is also briefly introduced, which will be
applied to tackle the optimization problems later.

\subsection{Design of Complementary Set of Sequences}

We are interested in developing efficient optimization methods for
the design of complementary sets of sequences. Consequently, to measure
the complementarity of a sequence set $\{\mathbf{x}_{m}\}_{m=1}^{M}$,
we consider the complementary integrated sidelobe level (CISL) metric
of a set of sequences, which is defined as 
\begin{equation}
{\rm CISL}=\sum_{k=1}^{N-1}\left|\sum_{m=1}^{M}r_{m,m}(k)\right|^{2}.\label{eq:ISL_set}
\end{equation}
Then a natural idea to generate complementary sets of unimodular sequences
is to minimize the CISL metric in \eqref{eq:ISL_set}, i.e., solving
the following optimization problem: 
\begin{equation}
\begin{array}{ll}
\underset{\{\mathbf{x}_{m}\}_{m=1}^{M}}{\mathsf{minimize}} & {\displaystyle \sum_{k=1}^{N-1}}\left|{\displaystyle \sum_{m=1}^{M}}r_{m,m}(k)\right|^{2}\\
\mathsf{subject\;to} & \left|x_{m}(n)\right|=1,\\
 & n=1,\ldots,N,\,m=1,\ldots,M.
\end{array}\label{eq:CSS_prob}
\end{equation}
Note that if the objective of problem \eqref{eq:CSS_prob} can be
driven to zero, then the corresponding solution is a complementary
set of sequences. But the problem may also be used to find almost
complementary sets of sequences for $(N,M)$ values for which no CSS
exists.

\subsection{Design of Sequence Set with Good Auto- and Cross-correlation Properties}

To design sequence sets with both good auto- and cross-correlation
properties, we consider the goodness measure used in \cite{SequenceSet_He2009},
which is defined as 
\begin{equation}
\Psi={\displaystyle \sum_{m=1}^{M}}{\displaystyle \sum_{\substack{k=1-N\\
k\neq0
}
}^{N-1}}\left|r_{m,m}(k)\right|^{2}+{\displaystyle \sum_{i=1}^{M}}{\displaystyle \sum_{\substack{j=1\\
j\neq i
}
}^{M}}{\displaystyle \sum_{k=1-N}^{N-1}}\left|r_{i,j}(k)\right|^{2}.\label{eq:obj_auto_cross}
\end{equation}
In this criterion, the first term contains the autocorrelation sidelobes
of all the sequences and the cross-correlations are involved in the
second term. Then, to design unimodular sequence sets with good correlation
properties, we consider the following optimization problem:
\begin{equation}
\begin{array}{ll}
\underset{\{\mathbf{x}_{m}\}_{m=1}^{M}}{\mathsf{minimize}} & {\displaystyle \sum_{m=1}^{M}}{\displaystyle \sum_{\substack{k=1-N\\
k\neq0
}
}^{N-1}}\left|r_{m,m}(k)\right|^{2}+{\displaystyle \sum_{i=1}^{M}}{\displaystyle \sum_{\substack{j=1\\
j\neq i
}
}^{M}}{\displaystyle \sum_{k=1-N}^{N-1}}\left|r_{i,j}(k)\right|^{2}\\
\mathsf{subject\;to} & \left|x_{m}(n)\right|=1,\,n=1,\ldots,N,\,m=1,\ldots,M.
\end{array}\label{eq:seqset_corr}
\end{equation}
Since $r_{m,m}(0)=N$, $m=1,\ldots,M$, due to the unimodular constraints,
problem \eqref{eq:seqset_corr} can be written more compactly as 
\begin{equation}
\begin{array}{ll}
\underset{\{\mathbf{x}_{m}\}_{m=1}^{M}}{\mathsf{minimize}} & {\displaystyle \sum_{i=1}^{M}}{\displaystyle \sum_{\substack{j=1}
}^{M}}{\displaystyle \sum_{k=1-N}^{N-1}}\left|r_{i,j}(k)\right|^{2}-N^{2}M\\
\mathsf{subject\;to} & \left|x_{m}(n)\right|=1,\\
 & n=1,\ldots,N,\,m=1,\ldots,M.
\end{array}\label{eq:seqset_prob}
\end{equation}
As have been shown in \cite{he2012waveform}, the criterion $\Psi$
defined in \eqref{eq:obj_auto_cross} is lower bounded by $N^{2}M(M-1)$
and thus cannot be made very small. This unveils the fact that it
is not possible to design a set of sequences with all auto- and cross-correlation
sidelobes very small. Therefore, we also consider the following more
general weighted formulation:
\begin{equation}
\begin{array}{ll}
\underset{\{\mathbf{x}_{m}\}_{m=1}^{M}}{\mathsf{minimize}} & {\displaystyle \sum_{i=1}^{M}}{\displaystyle \sum_{\substack{j=1}
}^{M}}{\displaystyle \sum_{k=1-N}^{N-1}}w_{k}\left|r_{i,j}(k)\right|^{2}-w_{0}N^{2}M\\
\mathsf{subject\;to} & \left|x_{m}(n)\right|=1,\,n=1,\ldots,N,\,m=1,\ldots,M,
\end{array}\label{eq:seqset_prob_weight}
\end{equation}
where $w_{k}=w_{-k}\geq0$, $k=0,\ldots,N-1$ are nonnegative weights
assigned to different time lags. It is easy to see that if we choose
$w_{k}=1$ for all $k,$ then problem \eqref{eq:seqset_prob_weight}
reduces to \eqref{eq:seqset_prob}. But problem \eqref{eq:seqset_prob_weight}
provides more flexibility in the sense that we can assign different
weights to different correlation lags, so that we can minimize the
correlations only within a certain time lag interval. Also note that
when $M=1$, problem \eqref{eq:seqset_prob_weight} becomes the weighted
integrated sidelobe level minimization problem considered in \cite{WISL_song2015}.

Two algorithms named CAN and WeCAN were proposed in \cite{SequenceSet_He2009}
to tackle problems \eqref{eq:seqset_prob_weight} and \eqref{eq:seqset_prob},
respectively. But the authors of \cite{SequenceSet_He2009} resorted
to solving ``almost equivalent'' problems that seem to work well
in practice. In this paper, we develop algorithms to directly tackle
the sequence set design formulations in \eqref{eq:seqset_prob_weight}
and \eqref{eq:seqset_prob}.

\subsection{The MM Method\label{sub:MM-Method}}

The MM method refers to the majorization-minimization method, which
is an approach to solve optimization problems that are too difficult
to solve directly. The principle behind the MM method is to transform
a difficult problem into a series of simple problems. Interested readers
may refer to \cite{hunter2004MMtutorial,MM_Stoica,razaviyayn2013unified}
and references therein for more details.

Suppose we want to minimize $f(\mathbf{x})$ over $\mathcal{X}\subseteq\mathbb{C}^{n}$.
Instead of minimizing the cost function $f(\mathbf{x})$ directly,
the MM approach optimizes a sequence of approximate objective functions
that majorize $f(\mathbf{x})$. More specifically, starting from a
feasible point $\mathbf{x}^{(0)},$ the algorithm produces a sequence
$\{\mathbf{x}^{(k)}\}$ according to the following update rule: 
\begin{equation}
\mathbf{x}^{(k+1)}\in\underset{\mathbf{x}\in\mathcal{X}}{\arg\min}\,\,u(\mathbf{x},\mathbf{x}^{(k)}),\label{eq:major_update}
\end{equation}
where $\mathbf{x}^{(k)}$ is the point generated by the algorithm
at iteration $k,$ and $u(\mathbf{x},\mathbf{x}^{(k)})$ is the majorization
function of $f(\mathbf{x})$ at $\mathbf{x}^{(k)}$. Formally, the
function $u(\mathbf{x},\mathbf{x}^{(k)})$ is said to majorize the
function $f(\mathbf{x})$ at the point $\mathbf{x}^{(k)}$ if 
\begin{eqnarray}
u(\mathbf{x},\mathbf{x}^{(k)}) & \geq & f(\mathbf{x}),\quad\forall\mathbf{x}\in\mathcal{X},\label{eq:major1}\\
u(\mathbf{x}^{(k)},\mathbf{x}^{(k)}) & = & f(\mathbf{x}^{(k)}).\label{eq:major2}
\end{eqnarray}
In other words, function $u(\mathbf{x},\mathbf{x}^{(k)})$ is an upper
bound of $f(\mathbf{x})$ over $\mathcal{X}$ and coincides with $f(\mathbf{x})$
at $\mathbf{x}^{(k)}$.

It is easy to show that with this scheme, the objective value is monotonically
decreasing (nonincreasing) at every iteration, i.e., 
\begin{equation}
f(\mathbf{x}^{(k+1)})\leq u(\mathbf{x}^{(k+1)},\mathbf{x}^{(k)})\leq u(\mathbf{x}^{(k)},\mathbf{x}^{(k)})=f(\mathbf{x}^{(k)}).\label{eq:descent-property}
\end{equation}
The first inequality and the third equality follow from the the properties
of the majorization function, namely \eqref{eq:major1} and \eqref{eq:major2}
respectively and the second inequality follows from \eqref{eq:major_update}.

To derive MM algorithms in practice, the key step is to find a majorization
function of the objective such that the majorized problem is easy
to solve. For that purpose, the following result on quadratic upper-bounding
will be useful later when constructing simple majorization functions.
\begin{lem}[\cite{MISL}]
\label{lem:majorizer}Let $\mathbf{L}$ be an $n\times n$ Hermitian
matrix and $\mathbf{M}$ be another $n\times n$ Hermitian matrix
such that $\mathbf{M}\succeq\mathbf{L}.$\textup{ }Then for any point
$\mathbf{x}_{0}\in\mathbf{C}^{n}$, the quadratic function $\mathbf{x}^{H}\mathbf{L}\mathbf{x}$
is majorized by $\mathbf{x}^{H}\mathbf{M}\mathbf{x}+2{\rm Re}\left(\mathbf{x}^{H}(\mathbf{L}-\mathbf{M})\mathbf{x}_{0}\right)+\mathbf{x}_{0}^{H}(\mathbf{M}-\mathbf{L})\mathbf{x}_{0}$
at\textup{ $\mathbf{x}_{0}$.}
\end{lem}

\section{Design of Complementary Set of Sequences via MM\label{sec:Seq-Design-viaMM}}

To tackle problem \eqref{eq:CSS_prob} via majorization-minimization,
we first perform some reformulations. Let us define an auxiliary
sequence of length $M(2N-1)$ as follows \cite{SequenceSet_Mojtaba2013}:
\begin{equation}
\mathbf{z}=[\mathbf{x}_{1}^{T},\mathbf{0}_{N-1}^{T},\ldots,\mathbf{x}_{M}^{T},\mathbf{0}_{N-1}^{T}]^{T},
\end{equation}
then the first $N$ aperiodic autocorrelation lags of $\mathbf{z}$
(denoted by $\{r_{z}(k)\}$) can be written as 
\begin{equation}
r_{z}(k)=\sum_{m=1}^{M}r_{m,m}(k),\,0\leq k\leq N-1.
\end{equation}
Then the sequence set $\{\mathbf{x}_{m}\}_{m=1}^{M}$ is complementary
if and only if $\mathbf{z}$ has a zero correlation zone (ZCZ) for
lags in the interval $1\leq k\leq N-1$, and the CSS design problem
\eqref{eq:CSS_prob} can be reformulated as 
\begin{equation}
\begin{array}{ll}
\underset{\{\mathbf{x}_{m}\}_{m=1}^{M}}{\mathsf{minimize}} & {\displaystyle \sum_{k=1}^{N-1}}\left|r_{z}(k)\right|^{2}\\
\mathsf{subject\;to} & \mathbf{z}=[\mathbf{x}_{1}^{T},\mathbf{0}_{N-1}^{T},\ldots,\mathbf{x}_{M}^{T},\mathbf{0}_{N-1}^{T}]^{T},\\
 & \left|x_{m}(n)\right|=1,\,n=1,\ldots,N,\,m=1,\ldots,M.
\end{array}\label{eq:CSS_WISL_prob}
\end{equation}
The objective in \eqref{eq:CSS_WISL_prob} can be viewed as the weighted
ISL metric in \cite{WISL_song2015} of the sequence $\mathbf{z}$
(i.e., $\sum_{k=1}^{M(2N-1)-1}w_{k}\left|r_{z}(k)\right|^{2}$) with
weights chosen as
\begin{equation}
w_{k}=\begin{cases}
1, & 1\leq k\leq N-1\\
0, & N\leq k\leq M(2N-1)-1.
\end{cases}\label{eq:weights}
\end{equation}
However, in problem \eqref{eq:CSS_WISL_prob}, the sequence $\mathbf{z}$
has some special structures and the original weighted ISL minimization
algorithm proposed in \cite{WISL_song2015} for designing unimodular
sequences cannot be directly applied due to the zeros. But the algorithm
can be adapted to take the sequence structure into account and in
the following we give a brief derivation of the modified algorithm,
which mainly follows from Section III.B in \cite{WISL_song2015}.

Similar to Section III.B in \cite{WISL_song2015}, we perform two
successive majorization steps to problem \eqref{eq:CSS_WISL_prob}.
Let $L=M(2N-1)$ be the length of $\mathbf{z}$, and $\mathbf{U}_{k},\,k=1-L,\ldots,,L-1$
be $L\times L$ Toeplitz matrices with the $k$th diagonal elements
being $1$ and $0$ elsewhere, i.e.,
\begin{equation}
\left[\mathbf{U}_{k}\right]_{i,j}=\begin{cases}
1 & \textrm{if }j-i=k\\
0 & \textrm{if }j-i\neq k,
\end{cases}\quad i,j=1,\ldots,L.\label{eq:U_k}
\end{equation}
Then the autocorrelations $\{r_{z}(k)\}$ of $\mathbf{z}$ can be
written in terms of $\mathbf{U}_{k}$ as 
\begin{equation}
r_{z}(k)=\mathbf{z}^{H}\mathbf{U}_{k}\mathbf{z},\,k=1-L,\ldots,,L-1.
\end{equation}
Then given $\mathbf{z}^{(l)}=[\mathbf{x}_{1}^{(l)T},\mathbf{0}_{N-1}^{T},\ldots,\mathbf{x}_{M}^{(l)T},\mathbf{0}_{N-1}^{T}]^{T}$
at iteration $l$, by using Lemma \ref{lem:majorizer} we can majorize
the objective of \eqref{eq:CSS_WISL_prob} by a quadratic function
as in \cite{WISL_song2015} and the majorized problem after the first
majorization step is given by

\begin{equation}
\begin{array}{ll}
\underset{\{\mathbf{x}_{m}\}_{m=1}^{M}}{\mathsf{minimize}} & \mathbf{z}^{H}\left(\mathbf{R}-(L-1)\mathbf{z}^{(l)}(\mathbf{z}^{(l)})^{H}\right)\mathbf{z}\\
\mathsf{subject\;to} & \mathbf{z}=[\mathbf{x}_{1}^{T},\mathbf{0}_{N-1}^{T},\ldots,\mathbf{x}_{M}^{T},\mathbf{0}_{N-1}^{T}]^{T},\\
 & \left|x_{m}(n)\right|=1,\,n=1,\ldots,N,\,m=1,\ldots,M,
\end{array}\label{eq:CSS_major1_prob}
\end{equation}
where 
\begin{eqnarray}
\mathbf{R} & = & {\displaystyle \sum_{\substack{k=1-L\\
k\neq0
}
}^{L-1}}w_{k}r_{z}^{(l)}(-k)\mathbf{U}_{k}
\end{eqnarray}
is a Hermitian Toeplitz matrix and $w_{k}=w_{-k}$, $k=1,\ldots,L-1$
are given in \eqref{eq:weights}. 

To perform the second majorization step, we first bound the maximum
eigenvalue of the matrix $\mathbf{R}-(L-1)\mathbf{z}^{(l)}(\mathbf{z}^{(l)})^{H}$
as in \cite{WISL_song2015}, i.e., 
\begin{equation}
\lambda_{{\rm max}}\left(\mathbf{R}-(L-1)\mathbf{z}^{(l)}(\mathbf{z}^{(l)})^{H}\right)\le\lambda_{u},\label{eq:eig_bound_L}
\end{equation}
where 
\begin{eqnarray}
\lambda_{u} & = & \frac{1}{2}\left(\max_{1\leq i\leq L}\mu_{2i}+\max_{1\leq i\leq L}\mu_{2i-1}\right),\\
\boldsymbol{\mu} & = & \mathbf{F}\mathbf{c},\label{eq:F_mu}\\
\mathbf{c} & = & [0,w_{1}r_{z}^{(l)}(1),\ldots,w_{L-1}r_{z}^{(l)}(L-1),\nonumber \\
 &  & 0,w_{L-1}r_{z}^{(l)}(1-L),\ldots,w_{1}r_{z}^{(l)}(-1)]^{T},
\end{eqnarray}
and the matrix $\mathbf{F}$ in \eqref{eq:F_mu} is the $2L\times2L$
FFT matrix with $F_{m,n}=e^{-j\frac{2mn\pi}{2L}},0\leq m,n<2L$. Then
by applying Lemma \ref{lem:majorizer} with $\mathbf{M=}\lambda_{u}\mathbf{I},$
we can obtain the majorized problem of \eqref{eq:CSS_major1_prob}
given by 
\begin{equation}
\begin{array}{ll}
\underset{\{\mathbf{x}_{m}\}_{m=1}^{M}}{\mathsf{minimize}} & {\rm Re}\left(\mathbf{z}^{H}\big(\mathbf{R}\!-\!(L-1)\mathbf{z}^{(l)}(\mathbf{z}^{(l)})^{H}\!-\!\lambda_{u}\mathbf{I}\big)\mathbf{z}^{(l)}\right)\\
\mathsf{subject\;to} & \mathbf{z}=[\mathbf{x}_{1}^{T},\mathbf{0}_{N-1}^{T},\ldots,\mathbf{x}_{M}^{T},\mathbf{0}_{N-1}^{T}]^{T},\\
 & \left|x_{m}(n)\right|=1,\,n=1,\ldots,N,\,m=1,\ldots,M,
\end{array}
\end{equation}
which can be rewritten as 
\begin{equation}
\begin{array}{ll}
\underset{\{\mathbf{x}_{m}\}_{m=1}^{M}}{\mathsf{minimize}} & \left\Vert \mathbf{z}-\mathbf{y}\right\Vert _{2}^{2}\\
\mathsf{subject\;to} & \mathbf{z}=[\mathbf{x}_{1}^{T},\mathbf{0}_{N-1}^{T},\ldots,\mathbf{x}_{M}^{T},\mathbf{0}_{N-1}^{T}]^{T},\\
 & \left|x_{m}(n)\right|=1,\,n=1,\ldots,N,\,m=1,\ldots,M,
\end{array}\label{eq:CSS_major2_prob}
\end{equation}
where 
\begin{eqnarray}
\mathbf{y} & = & -\big(\mathbf{R}-(L-1)\mathbf{z}^{(l)}(\mathbf{z}^{(l)})^{H}-\lambda_{u}\mathbf{I}\big)\mathbf{z}^{(l)}\nonumber \\
 & = & \left((L-1)MN+\lambda_{u}\right)\mathbf{z}^{(l)}-\mathbf{R}\mathbf{z}^{(l)}.
\end{eqnarray}
Problem \eqref{eq:CSS_major2_prob} admits the following closed form
solution
\begin{eqnarray}
x_{m}(n) & = & e^{j{\rm arg}(y_{(m-1)(2N-1)+n})},\nonumber \\
 &  & n=1,\ldots,N,m=1,\ldots,M.
\end{eqnarray}
The overall algorithm for the CSS design problem \eqref{eq:CSS_prob}
is summarized in Algorithm \ref{alg:CSS-MM}. Note that the algorithm
can be implemented by means of FFT (IFFT) operations, since $\mathbf{R}$
is Hermitian Toeplitz and it can be decomposed as 
\begin{equation}
\mathbf{R}=\frac{1}{2L}\mathbf{F}_{:,1:L}^{H}{\rm Diag}(\boldsymbol{\mu})\mathbf{F}_{:,1:L},
\end{equation}
according to Lemma 4 in \cite{WISL_song2015}.

\begin{algorithm}[tbh]
\begin{algor}[1]
\item [{Require:}] \begin{raggedright}
number of sequences $M$, sequence length $N$
\par\end{raggedright}
\item [{{*}}] \begin{raggedright}
Set $l=0$ and initialize $\{\mathbf{x}_{m}^{(0)}\}_{m=1}^{M}$. 
\par\end{raggedright}
\item [{{*}}] $L=M(2N-1)$
\item [{repeat}]~

\begin{algor}[1]
\item [{{*}}] $\mathbf{z}^{(l)}=[\mathbf{x}_{1}^{(l)T},\mathbf{0}_{N-1}^{T},\ldots,\mathbf{x}_{M}^{(l)T},\mathbf{0}_{N-1}^{T}]^{T}$
\item [{{*}}] $\mathbf{f}=\mathbf{F}[\mathbf{z}^{(l)T},\mathbf{0}_{1\times L}]^{T}$
\item [{{*}}] $\mathbf{r}=\frac{1}{2L}\mathbf{F}^{H}\left|\mathbf{f}\right|^{2}$
\item [{{*}}] $\mathbf{c}=\mathbf{r}\circ[0,\mathbf{1}_{N-1}^{T},\mathbf{0}_{2(L-N)+1}^{T},\mathbf{1}_{N-1}^{T}]^{T}$
\item [{{*}}] $\boldsymbol{\mu}=\mathbf{F}\mathbf{c}$
\item [{{*}}] $\lambda_{u}=\frac{1}{2}\big(\underset{1\leq i\leq N}{\max}\mu_{2i}+\underset{1\leq i\leq N}{\max}\mu_{2i-1}\big)$\ref{alg:MWISL-Set-adaptive}
\item [{{*}}] $\mathbf{y}=\left((L-1)MN+\lambda_{u}\right)\mathbf{z}^{(l)}-\frac{1}{2L}\mathbf{F}_{:,1:L}^{H}(\boldsymbol{\mu}\circ\mathbf{f})$
\item [{{*}}] $x_{m}^{(l+1)}(n)=e^{j{\rm arg}(y_{(m-1)(2N-1)+n})},n=1,\ldots,N,m=1,\ldots,M.$
\item [{{*}}] $l\leftarrow l+1$ 
\end{algor}
\item [{until}] convergence
\end{algor}
\protect\caption{\label{alg:CSS-MM}The MM Algorithm for CSS design problem \eqref{eq:CSS_prob}.}
\end{algorithm}

\section{Design of Sequence Set with Good Auto- and Cross-correlation Properties
via MM\label{sec:Set_Auto_Cross_Wei}}

In this section, we consider the problem of designing sequence sets
for both good auto- and cross-correlation properties. We first consider
the more general problem formulation with weights involved, i.e.,
problem \eqref{eq:seqset_prob_weight}, and derive an MM algorithm
for the problem in the following. 

Let us first stack the sequences $\mathbf{x}_{m},m=1,\ldots,M$ together
and denote it by $\mathbf{x}$, i.e., 
\begin{equation}
\mathbf{x}=[\mathbf{x}_{1}^{T},\ldots,\mathbf{x}_{M}^{T}]^{T},
\end{equation}
then we have 
\begin{equation}
\mathbf{x}_{m}=\mathbf{S}_{m}\mathbf{x},\,m=1,\ldots,M,\label{eq:x_m_x}
\end{equation}
where $\mathbf{S}_{m}$ is an $N\times NM$ block selection matrix
defined as 
\begin{equation}
\mathbf{S}_{m}=[\mathbf{0}_{N\times(m-1)N},\mathbf{I}_{N},\mathbf{0}_{N\times(M-m)N}].
\end{equation}
We then note that \eqref{eq:cross_corr-1} can be written more compactly
as 
\begin{equation}
r_{i,j}(k)=\mathbf{x}_{j}^{H}\mathbf{U}_{k}\mathbf{x}_{i},\,k=1-N,\ldots,N-1,\,i,j=1,\ldots,M,\label{eq:cross_corr_mat}
\end{equation}
where $\mathbf{U}_{k}$ is defined as in \eqref{eq:U_k} but is of
size $N\times N$ now. By combining \eqref{eq:cross_corr_mat} and
\eqref{eq:x_m_x}, we have 
\begin{eqnarray}
r_{i,j}(k) & = & \mathbf{x}^{H}\mathbf{S}_{j}^{H}\mathbf{U}_{k}\mathbf{S}_{i}\mathbf{x},\\
 &  & k=1-N,\ldots,N-1,\,i,j=1,\ldots,M,\nonumber 
\end{eqnarray}
and then
\begin{equation}
\begin{aligned}\left|r_{i,j}(k)\right|^{2} & =\left|\mathbf{x}^{H}\mathbf{S}_{j}^{H}\mathbf{U}_{k}\mathbf{S}_{i}\mathbf{x}\right|^{2}\\
 & =\left|{\rm Tr}\left(\mathbf{x}\mathbf{x}^{H}\mathbf{S}_{j}^{H}\mathbf{U}_{k}\mathbf{S}_{i}\right)\right|^{2}\\
 & =\left|{\rm vec}(\mathbf{x}\mathbf{x}^{H})^{H}{\rm vec}(\mathbf{S}_{j}^{H}\mathbf{U}_{k}\mathbf{S}_{i})\right|^{2}.
\end{aligned}
\label{eq:r_ij_square}
\end{equation}
By using \eqref{eq:r_ij_square}, problem \eqref{eq:seqset_prob_weight}
can be rewritten as 
\begin{equation}
\begin{array}{ll}
\underset{\mathbf{x}\in\mathbb{C}^{NM}}{\mathsf{minimize}} & {\rm vec}(\mathbf{x}\mathbf{x}^{H})^{H}\mathbf{L}{\rm vec}(\mathbf{x}\mathbf{x}^{H})-w_{0}N^{2}M\\
\mathsf{subject\;to} & \left|x_{n}\right|=1,\,n=1,\ldots,NM,
\end{array}\label{eq:prob_weight_vec}
\end{equation}
where 
\begin{equation}
\mathbf{L}={\displaystyle \sum_{i=1}^{M}}{\displaystyle \sum_{\substack{j=1}
}^{M}}{\displaystyle \sum_{k=1-N}^{N-1}}w_{k}{\rm vec}(\mathbf{S}_{j}^{H}\mathbf{U}_{k}\mathbf{S}_{i}){\rm vec}(\mathbf{S}_{j}^{H}\mathbf{U}_{k}\mathbf{S}_{i})^{H}.\label{eq:L_mat}
\end{equation}
Since $w_{k}\geq0,$ it is easy to see that $\mathbf{L}$ is a nonnegative
real symmetric matrix and it can be shown (see Lemma 5 in \cite{WISL_song2015})
that 
\begin{equation}
\mathbf{L}\preceq{\rm Diag}(\mathbf{b}),
\end{equation}
where $\mathbf{b}=\mathbf{L}\mathbf{1}$. Then given $\mathbf{x}^{(l)}$
at iteration $l$, by using Lemma \ref{lem:majorizer}, we know that
the objective of problem \eqref{eq:prob_weight_vec} is majorized
by the following function at $\mathbf{x}^{(l)}$: 
\begin{equation}
\begin{aligned} & \,\,u_{1}(\mathbf{x},\mathbf{x}^{(l)})\\
= & \,\,{\rm vec}(\mathbf{x}\mathbf{x}^{H})^{H}{\rm Diag}(\mathbf{b}){\rm vec}(\mathbf{x}\mathbf{x}^{H})\\
 & +2{\rm Re}\big({\rm vec}(\mathbf{x}\mathbf{x}^{H})^{H}(\mathbf{L}-{\rm Diag}(\mathbf{b})){\rm vec}(\mathbf{x}^{(l)}\mathbf{x}^{(l)H})\big)\\
 & +{\rm vec}(\mathbf{x}^{(l)}\mathbf{x}^{(l)H})^{H}({\rm Diag}(\mathbf{b})-\mathbf{L}){\rm vec}(\mathbf{x}^{(l)}\mathbf{x}^{(l)H})-w_{0}N^{2}M.
\end{aligned}
\label{eq:major_u1}
\end{equation}
Since the elements of $\mathbf{x}$ are of unit modulus, it is easy
to see that the first term of \eqref{eq:major_u1} is just a constant.
After ignoring the constant terms, the majorized problem of \eqref{eq:prob_weight_vec}
is given by 
\begin{equation}
\begin{array}{ll}
\underset{\mathbf{x}\in\mathbb{C}^{NM}}{\mathsf{minimize}} & {\rm Re}\big({\rm vec}(\mathbf{x}\mathbf{x}^{H})^{H}(\mathbf{L}-{\rm Diag}(\mathbf{b})){\rm vec}(\mathbf{x}^{(l)}\mathbf{x}^{(l)H})\big)\\
\mathsf{subject\;to} & \left|x_{n}\right|=1,\,n=1,\ldots,NM.
\end{array}\label{eq:prob_major1}
\end{equation}
By substituting $\mathbf{L}$ in \eqref{eq:L_mat} back, we have 
\begin{equation}
\begin{aligned} & {\rm Re}\big({\rm vec}(\mathbf{x}\mathbf{x}^{H})^{H}\mathbf{L}{\rm vec}(\mathbf{x}^{(l)}\mathbf{x}^{(l)H})\big)\\
= & {\displaystyle \sum_{i=1}^{M}}{\displaystyle \sum_{\substack{j=1}
}^{M}}{\displaystyle \sum_{k=1-N}^{N-1}}{\rm Re}\bigg(w_{k}{\rm Tr}\left(\mathbf{x}\mathbf{x}^{H}\mathbf{S}_{j}^{H}\mathbf{U}_{k}\mathbf{S}_{i}\right)\\
 & \qquad\qquad\qquad\quad\times{\rm Tr}\left(\mathbf{x}^{(l)}\mathbf{x}^{(l)H}\mathbf{S}_{i}^{H}\mathbf{U}_{-k}\mathbf{S}_{j}\right)\bigg)\\
= & {\displaystyle \sum_{i=1}^{M}}{\displaystyle \sum_{\substack{j=1}
}^{M}}{\displaystyle \sum_{k=1-N}^{N-1}}{\rm Re}\left(w_{k}r_{j,i}^{(l)}(-k)\mathbf{x}^{H}\mathbf{S}_{j}^{H}\mathbf{U}_{k}\mathbf{S}_{i}\mathbf{x}\right),
\end{aligned}
\label{eq:u1_term1}
\end{equation}
and the second term of the objective can also be rewritten as 
\begin{equation}
\begin{aligned} & {\rm Re}\big({\rm vec}(\mathbf{x}\mathbf{x}^{H})^{H}{\rm Diag}(\mathbf{b}){\rm vec}(\mathbf{x}^{(l)}\mathbf{x}^{(l)H})\big)\\
= & {\rm Re}\left({\rm vec}(\mathbf{x}\mathbf{x}^{H})^{H}\left(\mathbf{b}\circ{\rm vec}(\mathbf{x}^{(l)}\mathbf{x}^{(l)H})\right)\right)\\
= & {\rm Re}\left({\rm Tr}\left(\mathbf{x}\mathbf{x}^{H}{\rm mat}\left(\mathbf{b}\circ{\rm vec}(\mathbf{x}^{(l)}\mathbf{x}^{(l)H})\right)\right)\right)\\
= & {\rm Re}\left(\mathbf{x}^{H}\left({\rm mat}(\mathbf{b})\circ(\mathbf{x}^{(l)}\mathbf{x}^{(l)H})\right)\mathbf{x}\right),
\end{aligned}
\label{eq:u1_term2}
\end{equation}
where ${\rm mat}(\cdot)$ is the inverse operation of ${\rm vec}(\cdot)$.
It is clear that both \eqref{eq:u1_term1} and \eqref{eq:u1_term2}
are quadratic in $\mathbf{x}$ and problem \eqref{eq:prob_major1}
can be rewritten as
\begin{equation}
\begin{array}{ll}
\underset{\mathbf{x}\in\mathbb{C}^{NM}}{\mathsf{minimize}} & \mathbf{x}^{H}\left(\mathbf{R}-\mathbf{B}\circ\big(\mathbf{x}^{(l)}(\mathbf{x}^{(l)})^{H}\big)\right)\mathbf{x}\\
\mathsf{subject\;to} & \left|x_{n}\right|=1,\,n=1,\ldots,NM,
\end{array}\label{eq:prob_quad}
\end{equation}
where
\begin{equation}
\mathbf{R}={\displaystyle \sum_{i=1}^{M}}{\displaystyle \sum_{\substack{j=1}
}^{M}}{\displaystyle \sum_{k=1-N}^{N-1}}w_{k}r_{j,i}^{(l)}(-k)\mathbf{S}_{j}^{H}\mathbf{U}_{k}\mathbf{S}_{i},\label{eq:R_mat}
\end{equation}
\begin{equation}
\begin{aligned}\mathbf{B} & ={\rm mat}(\mathbf{b})\\
 & ={\rm mat}(\mathbf{L}\mathbf{1})\\
 & ={\rm mat}\left({\displaystyle \sum_{i=1}^{M}}{\displaystyle \sum_{\substack{j=1}
}^{M}}{\displaystyle \sum_{k=1-N}^{N-1}}w_{k}{\rm vec}(\mathbf{S}_{j}^{H}\mathbf{U}_{k}\mathbf{S}_{i}){\rm vec}(\mathbf{S}_{j}^{H}\mathbf{U}_{k}\mathbf{S}_{i})^{H}\mathbf{1}\right)\\
 & ={\rm mat}\left({\displaystyle \sum_{i=1}^{M}}{\displaystyle \sum_{\substack{j=1}
}^{M}}{\displaystyle \sum_{k=1-N}^{N-1}}w_{k}(N-\left|k\right|){\rm vec}(\mathbf{S}_{j}^{H}\mathbf{U}_{k}\mathbf{S}_{i})\right)\\
 & ={\displaystyle \sum_{i=1}^{M}}{\displaystyle \sum_{\substack{j=1}
}^{M}}{\displaystyle \sum_{k=1-N}^{N-1}}w_{k}(N-\left|k\right|)\mathbf{S}_{j}^{H}\mathbf{U}_{k}\mathbf{S}_{i}\\
 & =\mathbf{1}_{M\times M}\otimes\mathbf{W},
\end{aligned}
\label{eq:B_mat}
\end{equation}
and
\[
\begin{aligned}\mathbf{W} & ={\displaystyle \sum_{k=1-N}^{N-1}}w_{k}(N-\left|k\right|)\mathbf{U}_{k}\\
 & =\begin{bmatrix}w_{0}N & w_{1}(N-1) & \ldots & w_{N-1}\\
w_{1}(N-1) & w_{0}N & \ddots & \vdots\\
\vdots & \ddots & \ddots & w_{1}(N-1)\\
w_{N-1} & \ldots & w_{1}(N-1) & w_{0}N
\end{bmatrix}.
\end{aligned}
\]
 Note that in \eqref{eq:prob_quad} we have removed the ${\rm Re}(\cdot)$
operator since the matrices $\mathbf{R}$ and $\mathbf{B}$ are Hermitian.
Since the majorized problem \eqref{eq:prob_quad} is still hard to
solve directly, we propose to majorize the objective function at $\mathbf{x}^{(l)}$
again to further simplify the problem that we need to solve at each
iteration. Similarly, to construct a majorization function of the
quadratic objective in \eqref{eq:prob_quad}, we need to find a matrix
$\mathbf{M}$ such that $\mathbf{M}\succeq\mathbf{R}-\mathbf{B}\circ\big(\mathbf{x}^{(l)}(\mathbf{x}^{(l)})^{H}\big)$
and a straightforward choice may be $\mathbf{M}=\lambda_{{\rm max}}\left(\mathbf{R}-\mathbf{B}\circ\big(\mathbf{x}^{(l)}(\mathbf{x}^{(l)})^{H}\big)\right)\mathbf{I}$.
But to compute the maximum eigenvalue, some iterative algorithms are
needed and since we need to compute this at every iteration, it will
be computationally expensive. To maintain the computational efficiency
of the algorithm, here we propose to use some upper bound of $\lambda_{{\rm max}}\left(\mathbf{R}-\mathbf{B}\circ\big(\mathbf{x}^{(l)}(\mathbf{x}^{(l)})^{H}\big)\right)$
that can be easily computed. To derive such an upper bound, we first
introduce several results that will be useful. The first result reveals
a fact regarding the eigenvalues of the matrix $\mathbf{B}\circ\big(\mathbf{x}^{(l)}(\mathbf{x}^{(l)})^{H}\big)$,
which follows from \cite{WISL_song2015}.
\begin{lem}
\label{lem:eig_set}Let $\mathbf{B}$ be an $N\times N$ matrix and
$\mathbf{x}\in\mathbb{C}^{N}$ with $\left|x_{n}\right|=1,\,n=1,\ldots,N$.
Then $\mathbf{B}\circ(\mathbf{x}\mathbf{x}^{H})$ and $\mathbf{B}$
share the same set of eigenvalues. 
\end{lem}
The second result indicates some relations between the eigenvalues
of the Kronecker product of two matrices and the eigenvalues of the
two individual matrices \cite{roger1994topics}.
\begin{lem}
\label{lem:eig_kron}Let $\mathbf{A}$ and $\mathbf{B}$ be square
matrices of size $M$ and $N$, respectively. Let $\lambda_{1},\ldots,\lambda_{M}$
be the eigenvalues of $\mathbf{A}$ and $\mu_{1},\ldots,\mu_{N}$
be those of $\mathbf{B}$. Then the eigenvalues of $\mathbf{A}\otimes\mathbf{B}$
are $\lambda_{i}\mu_{j},i=1,\ldots,M$, $j=1,\ldots,N$ (including
algebraic multiplicities in all three cases).
\end{lem}
The third result regards bounds of the extreme eigenvalues of Hermitian
Toeplitz matrices, which can be computed by using FFTs \cite{eig_localization}. 
\begin{lem}
\label{lem:eig_bounds}Let $\mathbf{T}$ be an $N\times N$ Hermitian
Toeplitz matrix defined by $\{t_{k}\}_{k=0}^{N-1}$ as follows:
\[
\mathbf{T}=\begin{bmatrix}t_{0} & t_{1}^{*} & \ldots & t_{N-1}^{*}\\
t_{1} & t_{0} & \ddots & \vdots\\
\vdots & \ddots & \ddots & t_{1}^{*}\\
t_{N-1} & \ldots & t_{1} & t_{0}
\end{bmatrix}
\]
and $\mathbf{F}$ be a $2N\times2N$ FFT matrix with $F_{m,n}=e^{-j\frac{2mn\pi}{2N}},0\leq m,n<2N$.
Let $\mathbf{c}=[t_{0},t_{1},\cdots,t_{N-1},0,t_{N-1}^{*},\cdots,t_{1}^{*}]^{T}$
and $\boldsymbol{\mu}=\mathbf{F}\mathbf{c}$ be the discrete Fourier
transform of $\mathbf{c}$. Then 
\begin{eqnarray}
\lambda_{{\rm max}}(\mathbf{T}) & \leq & \frac{1}{2}\left(\max_{1\leq i\leq N}\mu_{2i}+\max_{1\leq i\leq N}\mu_{2i-1}\right),\label{eq:eig_up}\\
\lambda_{{\rm min}}(\mathbf{T}) & \geq & \frac{1}{2}\left(\min_{1\leq i\leq N}\mu_{2i}+\min_{1\leq i\leq N}\mu_{2i-1}\right).\label{eq:eig_low}
\end{eqnarray}

\end{lem}
Based on these results, we can now obtain an upper bound of $\lambda_{{\rm max}}\left(\mathbf{R}-\mathbf{B}\circ\big(\mathbf{x}^{(l)}(\mathbf{x}^{(l)})^{H}\big)\right)$
given in the following lemma.
\begin{lem}
\label{lem:eig_upperBound}Let $\mathbf{R}$ and $\mathbf{B}$ be
matrices defined in \eqref{eq:R_mat} and \eqref{eq:B_mat}, respectively.
Let $\mathbf{w}=[w_{0}N,w_{1}(N-1),\ldots,w_{N-1},0,w_{N-1},\ldots,w_{1}(N-1)]^{T},$
$\boldsymbol{\mu}=\mathbf{F}\mathbf{w}$ and $\lambda_{W}=\frac{1}{2}\left(\min_{1\leq i\leq N}\mu_{2i}+\min_{1\leq i\leq N}\mu_{2i-1}\right)$.
Then 
\begin{equation}
\lambda_{{\rm max}}\left(\mathbf{R}-\mathbf{B}\circ\big(\mathbf{x}^{(l)}(\mathbf{x}^{(l)})^{H}\big)\right)\leq\left\Vert \mathbf{R}\right\Vert -\lambda_{B},
\end{equation}
where 
\begin{equation}
\lambda_{B}=\begin{cases}
\min\left\{ M\lambda_{W},0\right\} , & M\geq2\\
\lambda_{W}, & M=1,
\end{cases}
\end{equation}
and $\left\Vert \cdot\right\Vert $ can be any submultiplicative matrix
norm.\end{lem}
\begin{IEEEproof}
See Appendix \ref{sec:Proof-of-Lemma-eigUB}.
\end{IEEEproof}
In our case, for computational efficiency, we choose the induced $\ell_{\infty}$-norm
(also known as max-row-sum norm) in Lemma \ref{lem:eig_upperBound},
which is defined as 
\begin{equation}
\left\Vert \mathbf{R}\right\Vert _{\infty}=\max_{i=1,\ldots,NM}\sum_{j=1}^{NM}\left|R_{i,j}\right|.
\end{equation}

Now, by choosing $\mathbf{M}=\left(\left\Vert \mathbf{R}\right\Vert _{\infty}-\lambda_{B}\right)\mathbf{I}$
in Lemma \ref{lem:majorizer}, the objective in \eqref{eq:prob_quad}
is majorized by 
\[
\begin{aligned} & \,\,u_{2}(\mathbf{x},\mathbf{x}^{(l)})\\
= & \,\,(\left\Vert \mathbf{R}\right\Vert _{\infty}-\lambda_{B})\mathbf{x}^{H}\mathbf{x}\\
 & +2{\rm Re}\left(\mathbf{x}^{H}\big(\mathbf{R}\!-\!\mathbf{B}\circ\big(\mathbf{x}^{(l)}(\mathbf{x}^{(l)})^{H}\big)\!-\!(\left\Vert \mathbf{R}\right\Vert _{\infty}\!-\!\lambda_{B})\mathbf{I}\big)\mathbf{x}^{(l)}\right)\\
 & +(\mathbf{x}^{(l)})^{H}((\left\Vert \mathbf{R}\right\Vert _{\infty}-\lambda_{B})\mathbf{I}-\mathbf{R}\!+\!\mathbf{B}\circ\big(\mathbf{x}^{(l)}(\mathbf{x}^{(l)})^{H}\big))\mathbf{x}^{(l)}.
\end{aligned}
\]
Again after ignoring the constant terms, the majorized problem of
\eqref{eq:prob_quad} is given by 
\begin{equation}
\begin{array}{ll}
\underset{\mathbf{x}\in\mathbb{C}^{NM}}{\mathsf{minimize}} & {\rm Re}\left(\mathbf{x}^{H}\mathbf{y}\right)\\
\mathsf{subject\;to} & \left|x_{n}\right|=1,\,n=1,\ldots,NM,
\end{array}\label{eq:prob_linear}
\end{equation}
where 
\begin{equation}
\mathbf{y}=\big(\mathbf{R}\!-\!\mathbf{B}\circ\big(\mathbf{x}^{(l)}(\mathbf{x}^{(l)})^{H}\big)\big)\mathbf{x}^{(l)}\!-\!(\left\Vert \mathbf{R}\right\Vert _{\infty}\!-\!\lambda_{B})\mathbf{x}^{(l)}.\label{eq:y}
\end{equation}
It is clear that problem \eqref{eq:prob_linear} is separable in the
elements of $\mathbf{x}$ and the solution of the problem is given
by 
\begin{equation}
x_{n}=e^{j{\rm arg}(-y_{n})},n=1,\ldots,NM.\label{eq:x_closed}
\end{equation}

According to the general steps of the majorization minimization method,
we can now implement the algorithm in a straightforward way, that
is at each iteration, we compute $\mathbf{y}$ according to \eqref{eq:y}
and update $\mathbf{x}$ via \eqref{eq:x_closed}. Clearly, the computational
cost is dominated by the computation of $\mathbf{y}$. To obtain an
efficient implementation, here we further explore the special structure
of the matrices involved in the computation of $\mathbf{y}$.

We first note that the matrix $\mathbf{R}$ in \eqref{eq:R_mat} can
be written as the following block matrix:
\begin{equation}
\mathbf{R}=\begin{bmatrix}\mathbf{R}_{11} & \mathbf{R}_{12} & \cdots & \mathbf{R}_{1M}\\
\mathbf{R}_{21} & \mathbf{R}_{22} & \cdots & \mathbf{R}_{2M}\\
\vdots & \vdots & \ddots & \vdots\\
\mathbf{R}_{M1} & \cdots & \cdots & \mathbf{R}_{MM}
\end{bmatrix},
\end{equation}
where each block is defined as
\begin{equation}
\mathbf{R}_{ij}={\displaystyle \sum_{k=1-N}^{N-1}}w_{k}r_{i,j}^{(l)}(-k)\mathbf{U}_{k},\,i,j=1,\ldots,M.
\end{equation}
It is easy to see that the building blocks $\mathbf{R}_{ij},i,j=1,\ldots,M$,
are Toeplitz matrices and when $i=j,$ they are also Hermitian. In
the following, we introduce a simple result regarding Toeplitz matrices
(not necessarily Hermitian) that can be used to perform the matrix
vector multiplication $\mathbf{R}\mathbf{x}^{(l)}$ more efficiently
via FFT (IFFT).
\begin{lem}
\label{lem:diagonal}Let $\mathbf{T}$ be an $N\times N$ Toeplitz
matrix defined as follows:
\[
\mathbf{T}=\begin{bmatrix}t_{0} & t_{1} & \ldots & t_{N-1}\\
t_{-1} & t_{0} & \ddots & \vdots\\
\vdots & \ddots & \ddots & t_{1}\\
t_{1-N} & \ldots & t_{-1} & t_{0}
\end{bmatrix}
\]
and $\mathbf{F}$ be a $2N\times2N$ FFT matrix with $F_{m,n}=e^{-j\frac{2mn\pi}{2N}},0\leq m,n<2N$.
Then $\mathbf{T}$ can be decomposed as $\mathbf{T}=\frac{1}{2N}\mathbf{F}_{:,1:N}^{H}{\rm Diag}(\mathbf{F}\mathbf{c})\mathbf{F}_{:,1:N}$,
where $\mathbf{c}=[t_{0},t_{-1},\cdots,t_{1-N},0,t_{N-1},\cdots,t_{1}]^{T}$.\end{lem}
\begin{IEEEproof}
See Appendix \ref{sec:Proof-of-Lemma-diag}.
\end{IEEEproof}
According to Lemma \ref{lem:diagonal}, by defining $\mathbf{H}$
to be the $2N\times N$ matrix composed of the first $N$ columns
of the $2N\times2N$ FFT matrix, i.e., 
\begin{equation}
\mathbf{H}=\mathbf{F}_{:,1:N},\label{eq:H_mat}
\end{equation}
we know that 
\begin{equation}
\mathbf{R}_{ij}=\frac{1}{2N}\mathbf{H}^{H}{\rm Diag}(\mathbf{F}\mathbf{c}_{ij})\mathbf{H},
\end{equation}
where
\begin{equation}
\begin{aligned}\mathbf{c}_{ij} & =[w_{0}r_{i,j}^{(l)}(0),w_{1}r_{i,j}^{(l)}(1),\ldots,w_{N-1}r_{i,j}^{(l)}(N-1),\\
 & \quad\,\,\,0,w_{N-1}r_{i,j}^{(l)}(1-N),\ldots,w_{1}r_{i,j}^{(l)}(-1)]^{T}.
\end{aligned}
\label{eq:c_ij}
\end{equation}
Thus, the matrix vector multiplication $\mathbf{R}\mathbf{x}^{(l)}$
can be performed as 
\begin{equation}
\mathbf{R}\mathbf{x}^{(l)}=\frac{1}{2N}\tilde{\mathbf{H}}^{H}\begin{bmatrix}{\rm Diag}(\mathbf{F}\mathbf{c}_{11}) & \cdots & {\rm Diag}(\mathbf{F}\mathbf{c}_{1M})\\
\vdots & \ddots & \vdots\\
{\rm Diag}(\mathbf{F}\mathbf{c}_{M1}) & \cdots & {\rm Diag}(\mathbf{F}\mathbf{c}_{MM})
\end{bmatrix}\tilde{\mathbf{H}}\mathbf{x}^{(l)},\label{eq:Rx_compute}
\end{equation}
where $\tilde{\mathbf{H}}$ is a $2MN\times MN$ block diagonal matrix
given by 
\begin{equation}
\tilde{\mathbf{H}}=\begin{bmatrix}\mathbf{H} & \mathbf{0} & \cdots & \mathbf{0}\\
\mathbf{0} & \mathbf{H} & \ddots & \vdots\\
\vdots & \ddots & \ddots & \mathbf{0}\\
\mathbf{0} & \cdots & \mathbf{0} & \mathbf{H}
\end{bmatrix}.
\end{equation}
From \eqref{eq:Rx_compute}, we can see that the multiplication $\mathbf{R}\mathbf{x}^{(l)}$
takes $M^{2}+2M$ FFT (IFFT) operations if all $\mathbf{c}_{ij},i,j=1,\ldots,M$
are given. Since to form the vectors $\mathbf{c}_{ij},i,j=1,\ldots,M,$
all the autocorrelations and cross-correlations, i.e., $r_{i,j}^{(l)}(k),i,j=1,\ldots,M,\,k=1-N,\ldots,N-1,$
are needed, and another $M^{2}$ FFT (IFFT) operations are required.
Similarly, $\left\Vert \mathbf{R}\right\Vert _{\infty}$ can also
be computed with $M^{2}+2M$ FFT (IFFT) operations, since it can be
obtained by taking the largest element of the vector $\tilde{\mathbf{R}}\mathbf{1}$,
where $\tilde{\mathbf{R}}$ is the matrix with each element being
the modulus of the corresponding element of $\mathbf{R}$, i.e., $\tilde{R}_{i,j}=\left|R_{i,j}\right|,i,j=1,\ldots,N.$
Finally, to compute $\big(\mathbf{B}\circ\big(\mathbf{x}^{(l)}(\mathbf{x}^{(l)})^{H}\big)\big)\mathbf{x}^{(l)}$
we first conduct some transformations as follows: 
\begin{equation}
\begin{aligned} & \,\,\big(\mathbf{B}\circ\big(\mathbf{x}^{(l)}(\mathbf{x}^{(l)})^{H}\big)\big)\mathbf{x}^{(l)}\\
= & \,\,{\rm diag}\big(\mathbf{B}{\rm Diag}(\mathbf{x}^{(l)})\big(\mathbf{x}^{(l)}(\mathbf{x}^{(l)})^{H}\big)^{T}\big)\\
= & \,\,{\rm diag}\big(\mathbf{B}\big(\mathbf{x}^{(l)}\circ\big(\mathbf{x}^{(l)}\big)^{*})\big(\mathbf{x}^{(l)}\big)^{T}\big)\\
= & \,\,{\rm diag}\big(\mathbf{B}\mathbf{1}_{NM\times1}\big(\mathbf{x}^{(l)}\big)^{T}\big)\\
= & \,\,(\mathbf{B}\mathbf{1}_{NM\times1})\circ\mathbf{x}^{(l)}\\
= & \,\,\left((\mathbf{1}_{M\times M}\otimes\mathbf{W})\mathbf{1}_{NM\times1}\right)\circ\mathbf{x}^{(l)}\\
= & \,\,\left(M\mathbf{1}_{M\times1}\otimes(\mathbf{W}\mathbf{1}_{N\times1})\right)\circ\mathbf{x}^{(l)}.
\end{aligned}
\end{equation}
Since $\mathbf{W}$ is Toeplitz, we know from Lemma \ref{lem:diagonal}
that it can be decomposed as 
\begin{equation}
\mathbf{W}=\frac{1}{2N}\mathbf{H}^{H}{\rm Diag}(\mathbf{F}\mathbf{w})\mathbf{H},
\end{equation}
where $\mathbf{w}$ is the same as the one defined in Lemma \ref{lem:eig_upperBound}.
Thus, $\big(\mathbf{B}\circ\big(\mathbf{x}^{(l)}(\mathbf{x}^{(l)})^{H}\big)\big)\mathbf{x}^{(l)}$
can be computed with $3$ FFT (IFFT) operations. 

In summary, to compute $\mathbf{y}$ as in \eqref{eq:y}, around $3M^{2}+4M+3$
$2N$-point FFT (IFFT) operations are needed. Since the computational
complexity of one FFT (IFFT) is $\mathcal{O}(N\log N),$ the per iteration
computational complexity of the proposed algorithm is of order $\mathcal{O}(M^{2}N\log N)$.
The overall algorithm is summarized in Algorithm \ref{alg:MWISL-Set}.

\begin{algorithm}[tbh]
\begin{algor}[1]
\item [{Require:}] \begin{raggedright}
number of sequences $M$, sequence length $N$, weights
$\{w_{k}\geq0\}_{k=0}^{N-1}$
\par\end{raggedright}
\item [{{*}}] Set $l=0$, initialize $\mathbf{x}^{(0)}$ of length $MN$. 
\item [{{*}}] \begin{raggedright}
$\mathbf{w}=[w_{0}N,w_{1}(N-1),\ldots,w_{N-1},0,w_{N-1},\ldots,w_{1}(N-1)]^{T}$
\par\end{raggedright}
\item [{{*}}] $\boldsymbol{\mu}=\mathbf{F}\mathbf{w}$
\item [{{*}}] \begin{raggedright}
$\lambda_{W}=\frac{1}{2}\left({\displaystyle \min_{1\leq i\leq N}}\mu_{2i}+{\displaystyle \min_{1\leq i\leq N}}\mu_{2i-1}\right)$
\par\end{raggedright}
\item [{{*}}] \begin{raggedright}
$\lambda_{B}=\begin{cases}
\min\left\{ M\lambda_{W},0\right\} , & M\geq2\\
\lambda_{W}, & M=1
\end{cases}$
\par\end{raggedright}
\item [{repeat}]~

\begin{algor}[1]
\item [{{*}}] \begin{raggedright}
Compute $r_{i,j}^{(l)}(k),i,j=1,\ldots,M,k=1-N,\dots,N-1$.
\par\end{raggedright}
\item [{{*}}] \begin{raggedright}
Compute $\mathbf{c}_{ij},i,j=1,\ldots,M$ according to \eqref{eq:c_ij}.
\par\end{raggedright}
\item [{{*}}] Compute $\mathbf{R}\mathbf{x}^{(l)}$ according to \eqref{eq:Rx_compute}.
\item [{{*}}] \begin{raggedright}
Compute $\left\Vert \mathbf{R}\right\Vert _{\infty}$ based
on $\left|\mathbf{c}_{ij}\right|,i,j=1,\ldots,M$.
\par\end{raggedright}
\item [{{*}}] $\mathbf{p}=\frac{M}{2N}\mathbf{1}_{M\times1}\otimes\left(\mathbf{H}^{H}\left(\boldsymbol{\mu}\circ(\mathbf{H}\mathbf{1})\right)\right)$
\item [{{*}}] $\mathbf{y}=\frac{\mathbf{R}\mathbf{x}^{(l)}-\mathbf{p}\circ\mathbf{x}^{(l)}}{\left\Vert \mathbf{R}\right\Vert _{\infty}\!-\!\lambda_{B}}\!-\!\mathbf{x}^{(l)}$
\item [{{*}}] $x_{n}^{(l+1)}=e^{j{\rm arg}(-y_{n})},\,n=1,\ldots,MN$ 
\item [{{*}}] $l\leftarrow l+1$ 
\end{algor}
\item [{until}] convergence
\end{algor}
\protect\caption{\label{alg:MWISL-Set}The MM Algorithm for problem \eqref{eq:seqset_prob_weight}.}
\end{algorithm}

\section{Simplified MM for the Case without Weights\label{sec:Adaptive-MM}}

In the previous section, we developed an algorithm for problem \eqref{eq:seqset_prob_weight}.
By simply choosing weights $w_{k}=1,k=1-N,\ldots,1+N$, the algorithm
can be readily applied to solve problem \eqref{eq:seqset_prob}. However,
as analyzed in the previous section, the algorithm requires about
$3M^{2}+4M$ $2N$-point FFT (IFFT) operations at every iteration.
In this section, we will derive an algorithm for problem \eqref{eq:seqset_prob},
which requires only $2M$ $2N$-point FFT (IFFT) operations per iteration.

Let us denote the sequence covariance matrix at lag $k$ by $\mathbf{R}_{k}$,
i.e., 
\begin{eqnarray}
\mathbf{R}_{k} & = & \begin{bmatrix}r_{1,1}(k) & r_{1,2}(k) & \ldots & r_{1,M}(k)\\
r_{2,1}(k) & r_{2,2}(k) &  & r_{2,M}(k)\\
\vdots &  & \ddots & \vdots\\
r_{M,1}(k) & \cdots & \cdots & r_{M,M}(k)
\end{bmatrix}\\
 &  & k=1-N,\ldots,,N-1.\nonumber 
\end{eqnarray}
By using \eqref{eq:cross_corr_mat}, it is easy to see that 
\begin{equation}
\mathbf{R}_{k}=\left(\mathbf{X}^{H}\mathbf{U}_{k}\mathbf{X}\right)^{T}=\mathbf{R}_{-k}^{H},\,k=0,\ldots,,N-1,
\end{equation}
where 
\begin{equation}
\mathbf{X}=[\mathbf{x}_{1},\ldots,\mathbf{x}_{M}].
\end{equation}
With the above matrix notation, problem \eqref{eq:seqset_prob} can
be rewritten as 
\begin{equation}
\begin{array}{ll}
\underset{\mathbf{X}\in\mathbb{C}^{N\times M}}{\mathsf{minimize}} & {\displaystyle \sum_{k=1-N}^{N-1}}\left\Vert \mathbf{X}^{H}\mathbf{U}_{k}\mathbf{X}\right\Vert _{F}^{2}-N^{2}M\\
\mathsf{subject\;to} & \left|X_{i,j}\right|=1,\,i=1,\ldots,N,\,j=1,\ldots,M.
\end{array}\label{eq:prob_set_mat}
\end{equation}
Since 
\begin{eqnarray*}
\left\Vert \mathbf{X}^{H}\mathbf{U}_{k}\mathbf{X}\right\Vert _{F}^{2} & = & {\rm Tr}\left(\mathbf{X}^{H}\mathbf{U}_{k}^{H}\mathbf{X}\mathbf{X}^{H}\mathbf{U}_{k}\mathbf{X}\right)\\
 & = & {\rm Tr}\left(\mathbf{X}\mathbf{X}^{H}\mathbf{U}_{k}^{H}\mathbf{X}\mathbf{X}^{H}\mathbf{U}_{k}\right)\\
 & = & {\rm vec}\left(\mathbf{X}\mathbf{X}^{H}\right)^{H}\left(\mathbf{U}_{k}^{H}\otimes\mathbf{U}_{k}^{H}\right){\rm vec}\left(\mathbf{X}\mathbf{X}^{H}\right),
\end{eqnarray*}
we have 
\begin{equation}
{\displaystyle \sum_{k=1-N}^{N-1}}\left\Vert \mathbf{X}^{H}\mathbf{U}_{k}\mathbf{X}\right\Vert _{F}^{2}={\rm vec}\left(\mathbf{X}\mathbf{X}^{H}\right)^{H}\tilde{\mathbf{L}}{\rm vec}\left(\mathbf{X}\mathbf{X}^{H}\right),\label{eq:obj_vec}
\end{equation}
where 
\begin{eqnarray}
\tilde{\mathbf{L}} & = & \sum_{k=1-N}^{N-1}\left(\mathbf{U}_{k}^{H}\otimes\mathbf{U}_{k}^{H}\right).\label{eq:L_tilde}
\end{eqnarray}
Let us define 
\begin{eqnarray}
\mathbf{h}_{p} & = & [1,e^{j\omega_{p}},\cdots,e^{j\omega_{p}(N-1)}]^{T},\,p=1,\ldots,2N,
\end{eqnarray}
where $\omega_{p}=\frac{2\pi}{2N}(p-1),\:p=1,\cdots,2N.$ Since $\mathbf{U}_{k}$
is Toeplitz and can be written in terms of $\mathbf{h}_{p},p=1,\ldots,2N$
according to Lemma \ref{lem:diagonal}, it can be shown that the matrix
$\tilde{\mathbf{L}}$ defined in \eqref{eq:L_tilde} can also be written
as 
\begin{equation}
\tilde{\mathbf{L}}=\frac{1}{2N}\sum_{p=1}^{2N}{\rm vec}(\mathbf{h}_{p}\mathbf{h}_{p}^{H}){\rm vec}(\mathbf{h}_{p}\mathbf{h}_{p}^{H})^{H},
\end{equation}
and then we have
\begin{equation}
\begin{aligned} & {\displaystyle \sum_{k=1-N}^{N-1}}\left\Vert \mathbf{X}^{H}\mathbf{U}_{k}\mathbf{X}\right\Vert _{F}^{2}\\
= & \frac{1}{2N}\sum_{p=1}^{2N}\left|{\rm vec}\left(\mathbf{X}\mathbf{X}^{H}\right)^{H}{\rm vec}(\mathbf{h}_{p}\mathbf{h}_{p}^{H})\right|^{2}\\
= & \frac{1}{2N}\sum_{p=1}^{2N}{\rm Tr}(\mathbf{X}\mathbf{X}^{H}\mathbf{h}_{p}\mathbf{h}_{p}^{H})^{2}\\
= & \frac{1}{2N}\sum_{p=1}^{2N}\left\Vert \mathbf{X}^{H}\mathbf{h}_{p}\right\Vert _{2}^{4}.
\end{aligned}
\end{equation}
Thus, problem \eqref{eq:prob_set_mat} can be further reformulated
as 
\begin{equation}
\begin{array}{ll}
\underset{\mathbf{X}\in\mathbb{C}^{N\times M}}{\mathsf{minimize}} & \frac{1}{2N}{\displaystyle \sum_{p=1}^{2N}}\left\Vert \mathbf{X}^{H}\mathbf{h}_{p}\right\Vert _{2}^{4}-N^{2}M\\
\mathsf{subject\;to} & \left|X_{i,j}\right|=1,\,i=1,\ldots,N,\,j=1,\ldots,M.
\end{array}\label{eq:prob_quartic}
\end{equation}

To construct a majorization function of the objective in \eqref{eq:prob_quartic},
we propose to majorize each$\left\Vert \mathbf{X}^{H}\mathbf{h}_{p}\right\Vert _{2}^{4}$
according to the following lemma.
\begin{lem}
\label{lem:quartic_major}Let $f(x)=x^{4},$ $x\in[0,t]$. Then for
given $x_{0}\in[0,t)$, $f(x)$ is majorized at $x_{0}$ over the
interval $[0,t]$ by the following quadratic function: 
\begin{equation}
ax^{2}+(4x_{0}^{3}-2ax_{0})x+ax_{0}^{2}-3x_{0}^{4},\label{eq:quad_major}
\end{equation}
where 
\begin{equation}
a=t^{2}+2x_{0}t+3x_{0}^{2}.\label{eq:a_val}
\end{equation}
\end{lem}
\begin{IEEEproof}
See Appendix \ref{sec:Proof-of-Lemma-quartic-major}.
\end{IEEEproof}
Given $\mathbf{X}^{(l)}$ at iteration $l$, by taking $\left\Vert \mathbf{X}^{H}\mathbf{h}_{p}\right\Vert _{2}$
as a whole, we know from Lemma \ref{lem:quartic_major} that each
$\left\Vert \mathbf{X}^{H}\mathbf{h}_{p}\right\Vert _{2}^{4}$ (for
any $p\in\{1,\ldots,2N\}$) is majorized by
\begin{equation}
a_{p}\left\Vert \mathbf{X}^{H}\mathbf{h}_{p}\right\Vert _{2}^{2}+b_{p}\left\Vert \mathbf{X}^{H}\mathbf{h}_{p}\right\Vert _{2}+a_{p}\left\Vert \mathbf{X}^{(l)H}\mathbf{h}_{p}\right\Vert _{2}^{2}-3\left\Vert \mathbf{X}^{(l)H}\mathbf{h}_{p}\right\Vert _{2}^{4},
\end{equation}
where 
\begin{eqnarray}
a_{p} & = & t^{2}+2t\left\Vert \mathbf{X}^{(l)H}\mathbf{h}_{p}\right\Vert _{2}+3\left\Vert \mathbf{X}^{(l)H}\mathbf{h}_{p}\right\Vert _{2}^{2},\label{eq:a_bound}\\
b_{p} & = & 4\left\Vert \mathbf{X}^{(l)H}\mathbf{h}_{p}\right\Vert _{2}^{3}-2a_{p}\left\Vert \mathbf{X}^{(l)H}\mathbf{h}_{p}\right\Vert _{2},
\end{eqnarray}
and $t$ is an upper bound of $\left\Vert \mathbf{X}^{H}\mathbf{h}_{p}\right\Vert _{2}^{4}$
over the set of interest at the current iteration. Since the objective
decreases at every iteration in the MM framework, at the current iteration
$l$, it is sufficient to consider the set on which the objective
is smaller than the current objective evaluated at $\mathbf{X}^{(l)}$.
Hence we can choose $t=\left({\displaystyle \sum_{p=1}^{2N}}\left\Vert \mathbf{X}^{(l)H}\mathbf{h}_{p}\right\Vert _{2}^{4}\right)^{1/4}$
here. Then the majorized problem of \eqref{eq:prob_quartic} is given
by (ignoring the constant terms and the scaling factor $\frac{1}{2N}$)
\begin{equation}
\begin{array}{ll}
\underset{\mathbf{X}\in\mathbb{C}^{N\times M}}{\mathsf{minimize}} & {\displaystyle \sum_{p=1}^{2N}}\left(a_{p}\left\Vert \mathbf{X}^{H}\mathbf{h}_{p}\right\Vert _{2}^{2}+b_{p}\left\Vert \mathbf{X}^{H}\mathbf{h}_{p}\right\Vert _{2}\right)\\
\mathsf{subject\;to} & \left|X_{i,j}\right|=1,\,i=1,\ldots,N,\,j=1,\ldots,M.
\end{array}\label{eq:quartic_major1}
\end{equation}

Let us first take a look at the first term of the objective. It can
be rewritten as follows: 
\begin{equation}
\begin{aligned}{\displaystyle \sum_{p=1}^{2N}}a_{p}\left\Vert \mathbf{X}^{H}\mathbf{h}_{p}\right\Vert _{2}^{2} & ={\displaystyle \sum_{p=1}^{2N}}a_{p}{\rm Tr}\left(\mathbf{X}^{H}\mathbf{h}_{p}\mathbf{h}_{p}^{H}\mathbf{X}\right)\\
 & ={\rm Tr}\left(\mathbf{X}^{H}\left(\sum_{p=1}^{2N}a_{p}\mathbf{h}_{p}\mathbf{h}_{p}^{H}\right)\mathbf{X}\right)\\
 & ={\rm Tr}\left(\mathbf{X}^{H}\mathbf{H}^{H}{\rm Diag}(\mathbf{a})\mathbf{H}\mathbf{X}\right),
\end{aligned}
\label{eq:first_term}
\end{equation}
where $\mathbf{H}=[\mathbf{h}_{1},\ldots,\mathbf{h}_{2N}]^{H}$ is
the matrix defined in \eqref{eq:H_mat} and $\mathbf{a}=[a_{1},\ldots,a_{2N}]^{T}$.
From Lemma \ref{lem:diagonal} and Lemma \ref{lem:eig_bounds}, we
can see that the matrix $\mathbf{H}^{H}{\rm Diag}(\mathbf{a})\mathbf{H}$
is Hermitian Toeplitz and its maximum eigenvalue is bounded above
as follows:
\begin{equation}
\lambda_{{\rm max}}(\mathbf{H}^{H}{\rm Diag}(\mathbf{a})\mathbf{H})\leq N\left(\max_{1\leq i\leq N}a_{2i}+\max_{1\leq i\leq N}a_{2i-1}\right).
\end{equation}
Let us define 
\begin{equation}
\lambda_{a}=N\left(\max_{1\leq i\leq N}a_{2i}+\max_{1\leq i\leq N}a_{2i-1}\right),
\end{equation}
then by choosing $\mathbf{M=}\lambda_{a}\mathbf{I}$ in Lemma \ref{lem:majorizer},
the function in \eqref{eq:first_term} is majorized by
\begin{equation}
\begin{aligned} & \,\,\lambda_{a}{\rm Tr}(\mathbf{X}^{H}\mathbf{X})\\
 & +2{\rm Re}\left({\rm Tr}\left(\mathbf{X}^{H}\left(\mathbf{H}^{H}{\rm Diag}(\mathbf{a})\mathbf{H}-\lambda_{a}\mathbf{I}\right)\mathbf{X}^{(l)}\right)\right)\\
 & +{\rm Tr}\left(\mathbf{X}^{(l)H}\left(\lambda_{a}\mathbf{I}-\mathbf{H}^{H}{\rm Diag}(\mathbf{a})\mathbf{H}\right)\mathbf{X}^{(l)}\right).
\end{aligned}
\label{eq:first_term_major}
\end{equation}
Note that ${\rm Tr}(\mathbf{X}^{H}\mathbf{X})=MN,$ so the first term
of \eqref{eq:first_term_major} is just a constant. 

For the second term of the objective in \eqref{eq:quartic_major1},
we have 
\begin{equation}
\begin{array}{cl}
 & {\displaystyle \sum_{p=1}^{2N}}b_{p}\left\Vert \mathbf{X}^{H}\mathbf{h}_{p}\right\Vert _{2}\\
= & {\displaystyle \sum_{p=1}^{2N}}\left(4\left\Vert \mathbf{X}^{(l)H}\mathbf{h}_{p}\right\Vert _{2}^{2}-2a_{p}\right)\left\Vert \mathbf{X}^{(l)H}\mathbf{h}_{p}\right\Vert _{2}\left\Vert \mathbf{X}^{H}\mathbf{h}_{p}\right\Vert _{2}\\
\text{\ensuremath{\leq}} & {\displaystyle \sum_{p=1}^{2N}}\left(4\left\Vert \mathbf{X}^{(l)H}\mathbf{h}_{p}\right\Vert _{2}^{2}-2a_{p}\right){\rm Re}\left(\mathbf{h}_{p}^{H}\mathbf{X}^{(l)}\mathbf{X}^{H}\mathbf{h}_{p}\right)\\
= & {\rm Re}\left({\rm Tr}\left(\tilde{\mathbf{Y}}\mathbf{X}^{H}\right)\right)
\end{array}\label{eq:bX_major}
\end{equation}
where 
\begin{equation}
\tilde{\mathbf{Y}}=\left({\displaystyle \sum_{p=1}^{2N}}\left(4\left\Vert \mathbf{X}^{(l)H}\mathbf{h}_{p}\right\Vert _{2}^{2}-2a_{p}\right)\mathbf{h}_{p}\mathbf{h}_{p}^{H}\right)\mathbf{X}^{(l)}
\end{equation}
and the inequality follows from the Cauchy-Schwarz inequality and
the fact 
\begin{equation}
\begin{aligned}4\left\Vert \mathbf{X}^{(l)H}\mathbf{h}_{p}\right\Vert _{2}^{2}-2a_{p} & =-2\left(\left\Vert \mathbf{X}^{(l)H}\mathbf{h}_{p}\right\Vert _{2}+t\right)^{2}\leq0.\end{aligned}
\end{equation}
Since the inequality in \eqref{eq:bX_major} holds with equality when
$\mathbf{X}=\mathbf{X}^{(l)}$, ${\rm Re}\left({\rm Tr}\left(\tilde{\mathbf{Y}}\mathbf{X}^{H}\right)\right)$
majorizes the second term of the objective in \eqref{eq:quartic_major1}
at $\mathbf{X}^{(l)}$.

By adding the two majorization functions, i.e., \eqref{eq:first_term_major}
and \eqref{eq:bX_major}, we get the majorized problem of \eqref{eq:quartic_major1}
(ignoring the constant terms): 

\begin{equation}
\begin{array}{ll}
\underset{\mathbf{X}\in\mathbb{C}^{N\times M}}{\mathsf{minimize}} & {\rm Re}\left({\rm Tr}\left(\mathbf{Y}\mathbf{X}^{H}\right)\right)\\
\mathsf{subject\;to} & \left|X_{i,j}\right|=1,\,i=1,\ldots,N,\,j=1,\ldots,M,
\end{array}\label{eq:prob_linear_X}
\end{equation}
where
\begin{equation}
\begin{aligned}\mathbf{Y} & =\tilde{\mathbf{Y}}+2\left(\mathbf{H}^{H}{\rm Diag}(\mathbf{a})\mathbf{H}-\lambda_{a}\mathbf{I}\right)\mathbf{X}^{(l)}\\
 & =4\left({\displaystyle \sum_{p=1}^{2N}}\left\Vert \mathbf{X}^{(l)H}\mathbf{h}_{p}\right\Vert _{2}^{2}\mathbf{h}_{p}\mathbf{h}_{p}^{H}\right)\mathbf{X}^{(l)}-2\lambda_{a}\mathbf{X}^{(l)}.
\end{aligned}
\label{eq:Y_mat}
\end{equation}
 It is easy to see that problem \eqref{eq:prob_linear_X} can be rewritten
as 
\begin{equation}
\begin{array}{ll}
\underset{\mathbf{X}\in\mathbb{C}^{N\times M}}{\mathsf{minimize}} & {\displaystyle \sum_{i=1}^{N}\sum_{j=1}^{M}{\rm Re}\left(X_{i,j}^{*}Y_{i,j}\right)}\\
\mathsf{subject\;to} & \left|X_{i,j}\right|=1,\,i=1,\ldots,N,\,j=1,\ldots,M,
\end{array}
\end{equation}
which is separable in the elements of $\mathbf{X}$ and the solution
of the problem is given by 
\begin{equation}
X_{i,j}=e^{j{\rm arg}(-Y_{i,j})},i=1,\ldots,N,\,j=1,\ldots,M.\label{eq:X_closed}
\end{equation}

Then at every iteration of the algorithm, we just compute the matrix
$\mathbf{Y}$ given in \eqref{eq:Y_mat} and update $\mathbf{X}$
according to \eqref{eq:X_closed}. It is worth noting that the matrix
$\mathbf{Y}$ in \eqref{eq:Y_mat} can be computed efficiently via
FFT (IFFT), since it can be rewritten as 
\begin{equation}
\mathbf{Y}=4\mathbf{H}^{H}{\rm Diag}(\mathbf{q})\mathbf{H}\mathbf{X}^{(l)}-2\lambda_{a}\mathbf{X}^{(l)},\label{eq:Y_mat_FFT}
\end{equation}
where 
\begin{equation}
\mathbf{q}=\left|\mathbf{H}\mathbf{X}^{(l)}\right|^{2}\mathbf{1}_{M\times1}
\end{equation}
and $\left|\cdot\right|^{2}$ denotes the element-wise absolute-squared
value. The overall algorithm is then summarized in Algorithm \ref{alg:MWISL-Set-adaptive}
and we can see that $2M$ 2N-point FFT (IFFT) operations are needed
at each iteration.

\begin{algorithm}[tbh]
\begin{algor}[1]
\item [{Require:}] \begin{raggedright}
number of sequences $M$, sequence length $N$
\par\end{raggedright}
\item [{{*}}] Set $l=0$, initialize $\mathbf{X}^{(0)}$ of size $N\times M$. 
\item [{repeat}]~

\begin{algor}[1]
\item [{{*}}] $\mathbf{q}=\left|\mathbf{H}\mathbf{X}^{(l)}\right|^{2}\mathbf{1}_{M\times1}$
\item [{{*}}] $t=\left(\mathbf{1}^{T}\left(\mathbf{q}\circ\mathbf{q}\right)\right)^{\frac{1}{4}}$
\item [{{*}}] \begin{raggedright}
$a_{i}=t^{2}+2t\sqrt{q_{i}}+3q_{i},i=1,\ldots,2N$
\par\end{raggedright}
\item [{{*}}] $\lambda_{a}=N\left({\displaystyle \max_{1\leq i\leq N}}a_{2i}+{\displaystyle \max_{1\leq i\leq N}}a_{2i-1}\right)$\smallskip{}

\item [{{*}}] $\mathbf{Y}=4\mathbf{H}^{H}{\rm Diag}(\mathbf{q})\mathbf{H}\mathbf{X}^{(l)}-2\lambda_{a}\mathbf{X}^{(l)}$
\item [{{*}}] \begin{raggedright}
$X_{i,j}^{(l+1)}\negthinspace=\negthinspace e^{j{\rm arg}(-Y_{i,j})},i=1,\ldots,N,j=1,\ldots,M$
\par\end{raggedright}
\item [{{*}}] $l\leftarrow l+1$ 
\end{algor}
\item [{until}] convergence
\end{algor}
\protect\caption{\label{alg:MWISL-Set-adaptive}The MM Algorithm for problem \eqref{eq:seqset_prob}.}
\end{algorithm}

\section{Convergence Analysis and Acceleration Scheme\label{sec:Convergence-Acc}}

\subsection{Convergence Analysis}

The algorithms developed in the previous sections are all based on
the general majorization-minimization method and according to subsection
\ref{sub:MM-Method} we know that the sequences of objective values
generated by the algorithms at every iteration are nonincreasing.
Since it is easy to see that the objective functions of problems \eqref{eq:CSS_prob},
\eqref{eq:seqset_prob} and \eqref{eq:seqset_prob_weight} are all
bounded below by $0,$ the sequences of objective values are guaranteed
to converge to finite values.

In the following, we establish the convergence of the solution sequences
generated by the algorithms to stationary points. Let $f(\mathbf{x})$
be a differentiable function and $\mathcal{X}$ be an arbitrary constraint
set, then a point $\mathbf{x}^{\star}\in\mathcal{X}$ is said to be
a stationary point of the problem
\begin{equation}
\begin{array}{ll}
\underset{\mathbf{x}\in\mathcal{X}}{\mathsf{minimize}} & f(\mathbf{x})\end{array}
\end{equation}
if it satisfies the following first-order optimality condition \cite{Bertsekas2003}:
\[
\nabla f(\mathbf{x}^{\star})^{T}\mathbf{z}\geq0,\,\forall\mathbf{z}\in T_{\mathcal{X}}(\mathbf{x}^{\star}),
\]
where $T_{\mathcal{X}}(\mathbf{x}^{\star})$ denotes the tangent cone
of $\mathcal{X}$ at $\mathbf{x}^{\star}.$ The convergence property
of the CSS design algorithm in Algorithm \ref{alg:CSS-MM} can be
stated as follows.
\begin{thm}
\label{thm:MM-converge}Let $\{\mathbf{x}_{m}^{(l)}\}_{m=1}^{M},l=0,1,\ldots$
be the sequence of iterates generated by Algorithm \ref{alg:CSS-MM}.
Then the sequence has at least one limit point and every limit point
of the sequence is a stationary point of problem \eqref{eq:CSS_prob}. \end{thm}
\begin{IEEEproof}
The proof is similar to that given in \cite{WISL_song2015} and we
omit it here.
\end{IEEEproof}
Note that the convergence results of Algorithms \ref{alg:MWISL-Set}
and \ref{alg:MWISL-Set-adaptive} can be stated similarly and the
sequences generated by the two algorithms converge to stationary points
of problems \eqref{eq:seqset_prob_weight} and \eqref{eq:seqset_prob},
respectively.

\subsection{Acceleration Scheme \label{sec:Acceleration-Schemes}}

The popularity of the MM method is due to its simplicity and numerical
stability (monotonicity), but it is usually attained at the expense
of slow convergence. Due to the successive majorization steps that
we have carried out in the derivation of the majorization functions,
the convergence of the proposed algorithms seems to be slow. To fix
this issue, we can apply some acceleration schemes and in this subsection
we briefly introduce such a scheme that can be easily applied to speed
up the proposed MM algorithms. It is the squared iterative method
(SQUAREM) \cite{SQUAREM}, which was originally proposed to accelerate
any Expectation\textendash Maximization (EM) algorithms. It seeks
to approximate Newton\textquoteright s method for finding a fixed
point of the EM algorithm map and generally achieves superlinear convergence.
Since SQUAREM only requires the EM updating map, it can be readily
applied to any EM-type algorithms. In \cite{WISL_song2015}, it was
applied to accelerate some MM algorithms and some modifications were
made to maintain the monotonicity of the original MM algorithm and
to ensure the feasibility of the solution after every iteration. The
modified scheme is summarized in Algorithm 3 in \cite{WISL_song2015}
and we will apply it to accelerate the proposed MM algorithms in this
paper.

\section{Numerical Experiments\label{sec:Numerical-Experiments}}

To show the performance of the proposed algorithms in designing set
of sequences for various scenarios, we present some experimental results
in this section. For clarity, the MM algorithms proposed for problems
\eqref{eq:CSS_prob}, \eqref{eq:seqset_prob} and \eqref{eq:seqset_prob_weight},
i.e., Algorithms \ref{alg:CSS-MM}, \ref{alg:MWISL-Set-adaptive}
and \ref{alg:MWISL-Set}, will be referred to as MM-CSS, MM-Corr and
MM-WeCorr, respectively. And the acceleration scheme described in
section \ref{sec:Acceleration-Schemes} was applied in our implementation
of the algorithms. All experiments were performed in Matlab on a PC
with a 3.20 GHz i5-3470 CPU and 8 GB RAM.

\subsection{CSS Design}

\begin{figure}[htbp]
\centering{}\includegraphics[width=0.95\columnwidth]{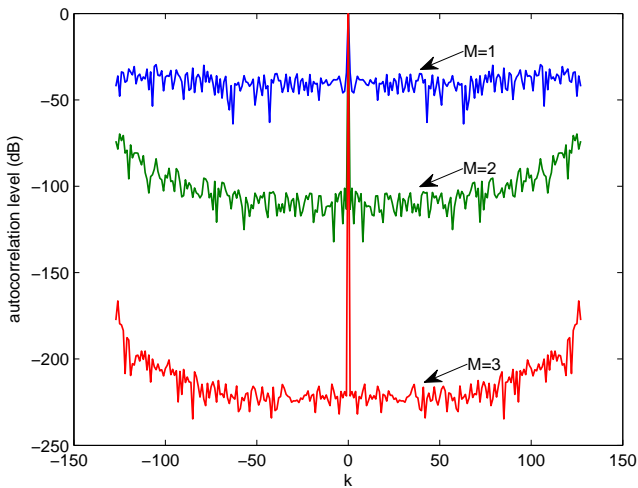}\protect\caption{\label{fig:CSS_corr}Autocorrelation levels of sequence sets with
$N=128$ and $M=1,2,3$.}
\end{figure}

In this subsection, we give an example of applying the proposed MM-CSS
algorithm to design (almost) complementary sets of sequences (CSS).
We consider the design of unimodular CSS of length $N=128$ and with
$M=1,2,3.$ For all cases, the initial sequence set $\{\mathbf{x}_{m}^{(0)}\}_{m=1}^{M}$
was generated randomly with each sequence being $\{e^{j2\pi\theta_{n}}\}_{n=1}^{N}$,
where $\{\theta_{n}\}_{n=1}^{N}$ are independent random variables
uniformly distributed in $[0,1]$. The stopping criterion was set
to be $\left|{\rm ISL}^{(l+1)}-{\rm ISL}^{(l)}\right|/\max\left(1,{\rm ISL}^{(l)}\right)\leq10^{-15}$
to allow enough iterations. The complementary autocorrelation levels
of the output sequence sets with $M=1,2,3$ sequences are shown in
Fig. \ref{fig:CSS_corr}, where the complementary autocorrelation
level is the normalized autocorrelation sum in dB defined as

\begin{equation}
20\log_{10}\frac{\left|\sum_{m=1}^{M}r_{m,m}(k)\right|}{\sum_{m=1}^{M}r_{m,m}(0)},\,k=1-N,\ldots,N-1.
\end{equation}
From the figure, we can see that as $M$ increases, the complementary
autocorrelation level decreases, which can be easily understood as
larger $M$ provides more degrees of freedom for the CSS design. In
particular, when $M=3$ the autocorrelation sums of the sequences
are very close to zero and the sequences can be viewed as complementary
in practice. 

\subsection{Approaching the Lower Bound of $\Psi$}

As have been mentioned earlier, the criterion $\Psi$ defined in \eqref{eq:obj_auto_cross}
is lower bounded by $N^{2}M(M-1)$. Then a natural question is whether
we can achieve that bound. In this subsection, we apply the proposed
MM-Corr and MM-WeCorr algorithms to minimize the criterion $\Psi$,
i.e., solving problem \eqref{eq:seqset_prob}, and compare the performance
with the CAN algorithm \cite{SequenceSet_He2009}.

In the experiment, we consider sequences sets with $M\in\{2,3,4\}$
sequences and each sequence of length $N\in\{256,1024\}$. For all
algorithms, the initial sequence set was generated randomly as in
the previous subsection, and the stopping criterion was set to be
$\left|\Psi^{(l+1)}-\Psi^{(l)}\right|/\Psi^{(l)}\leq10^{-8}$.

For each $(M,N)$ pair, the algorithms were repeated 10 times and
the minimum and average values of $\Psi$ achieved by the three algorithms,
together with the corresponding lower bound, are shown in Table \ref{tab:Corr_bound}.
The average running time of the three algorithms was also recorded
and is provided in Table \ref{tab:avg-running-time}. From Table
\ref{tab:Corr_bound}, we can see that all the three algorithms can
get reasonably close to the lower bound of $\Psi$, which means the
sequence sets generated by the algorithms are almost optimal for the
$(M,N)$ pairs that have been considered. Another point we notice
is that, for all $(M,N)$ pairs and all algorithms, the average values
over 10 random trials are quite close to the minimum values, which
implies that the three algorithms are not sensitive to the initial
points. From Table \ref{tab:avg-running-time}, we can see that for
each $(M,N)$ pair, the MM-Corr algorithm is the fastest and the CAN
algorithm is the slowest among the three algorithms. Since the per
iteration computational complexity of MM-Corr and CAN is almost the
same ($2M$ $2N$-point FFT (IFFT) operations), it implies that MM-Corr
takes far fewer iterations to converge compared with CAN. Another
observation is that for the same sequence length $N$, the cases with
larger $M$ values take less time compared with the cases with smaller
$M$ values, for example the running time of the algorithms for the
pair $(M=4,N=256)$ is less than that for $(M=2,N=256)$. Since a
larger $M$ value means higher per iteration computational complexity,
the observation implies that when $M$ becomes larger, the algorithms
need much fewer iterations to converge. It probably further implies
that it is easier for a larger set of sequences to approach the lower
bound than a smaller set of sequences.

\begin{table*}[t]
\protect\caption{\label{tab:Corr_bound}The lower bound of $\Psi$ in \eqref{eq:obj_auto_cross}
and the values achieved by different algorithms.}

\begin{centering}
\begin{tabular*}{0.9\textwidth}{@{\extracolsep{\fill}}lccccccc}
\hline 
 & \multicolumn{2}{c}{CAN} & \multicolumn{2}{c}{MM-WeCorr} & \multicolumn{2}{c}{MM-Corr} & Lower Bound\tabularnewline
\hline 
 & minimum & average & minimum & average & minimum & average & \tabularnewline
\hline 
$M=2,N=256$ & 131082 & 131089 & 131083 & 131093 & \textbf{131079} & 131093 & 131072\tabularnewline
$M=3,N=256$ & 393220 & 393222 & \textbf{393217} & 393220 & 393219 & 393222 & 393216\tabularnewline
$M=4,N=256$ & 786436 & 786439 & \textbf{786433} & 786436 & \textbf{786433} & 786436 & 786432\tabularnewline
$M=2,N=1024$ & 2097336 & 2097394 & 2097426 & 2098298 & \textbf{2097335} & 2097453 & 2097152\tabularnewline
$M=3,N=1024$ & 6291553 & 6291580 & \textbf{6291486} & 6291556 & 6291504 & 6291548 & 6291456\tabularnewline
$M=4,N=1024$ & 12582992 & 12583019 & \textbf{12582937} & 12582989 & 12582939 & 12582992 & 12582912\tabularnewline
\hline 
\end{tabular*}
\par\end{centering}

\end{table*}

\begin{table}[tbh]
\protect\caption{\label{tab:avg-running-time}The average running time (in seconds)
of different algorithms over 10 random trials.}

\centering{}%
\begin{tabular*}{0.95\columnwidth}{@{\extracolsep{\fill}}lccc}
\hline 
 & CAN & MM-WeCorr & MM-Corr\tabularnewline
\hline 
$M=2,N=256$ & 9.3342 & 0.6765 & 0.2435\tabularnewline
$M=3,N=256$ & 2.3461 & 0.3813 & 0.1000\tabularnewline
$M=4,N=256$ & 1.3562 & 0.3822 & 0.0844\tabularnewline
$M=2,N=1024$ & 33.8459 & 1.2011 & 0.6137\tabularnewline
$M=3,N=1024$ & 8.0584 & 1.0797 & 0.2750\tabularnewline
$M=4,N=1024$ & 4.9846 & 1.0298 & 0.2242\tabularnewline
\hline 
\end{tabular*}
\end{table}

\subsection{Sequence Set Design with Zero Correlation Zone }

As can be seen from the previous subsection, it is impossible to design
a set of sequences with all auto- and cross-correlation sidelobes
very small. Since in some applications, it is enough to minimize the
correlations only within a certain time lag interval, in this subsection
we present an example of applying the proposed MM-WeCorr algorithm
to design a set of sequences with low correlation sidelobes only at
required lags and compare the performance with the WeCAN algorithm
in \cite{SequenceSet_He2009}. The Matlab code of the WeCAN algorithm
was downloaded from the website\footnote{http://www.sal.ufl.edu/book/}
of the book \cite{he2012waveform}.

\begin{figure}[t]
\centering{}\includegraphics[width=0.95\columnwidth]{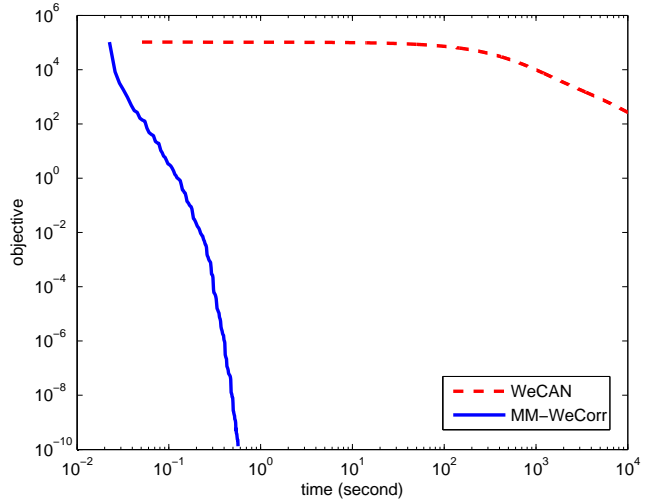}\protect\caption{\label{fig:obj_vs_time}Evolution of the objective with respect to
the running time (in seconds).}
\end{figure}

Suppose we want to design a sequence set with $M=3$ sequences each
of length $N=256$ and with low auto- and cross-correlations only
at lags $k=51,\ldots,80$. To tackle the problem, we apply the MM-WeCorr
and WeCAN algorithms from random initial sequence sets generated as
in the previous subsections. For the MM-WeCorr algorithm, we choose
the weights $\{w_{k}\}_{k=0}^{N-1}$ as follows:

\begin{equation}
w_{k}=\begin{cases}
1, & k\in\{51,\dots,80\}\\
0, & {\rm otherwise},
\end{cases}\label{eq:weights-sim}
\end{equation}
so that only the correlations at the required lags will be minimized.
For both algorithms, we do not stop until the objective in \eqref{eq:seqset_prob_weight}
goes below $10^{-10}$ or after 10000 seconds. The evolution curves
of the objective with respect to the running time are shown in Fig.
\ref{fig:obj_vs_time}. From the figure we can see that the proposed
MM-WeCorr algorithm drives the objective to $10^{-10}$ within $1$
second, while the objective is still above $10^{2}$ after 10000 seconds
for WeCAN. This is because the proposed MM-WeCorr algorithm requires
about $3M^{2}+4M$ $2N$-point FFT's per iteration, while each iteration
of WeCAN requires $2MN$ computations of $2N$-point FFT's and also
$2N$ computations of the SVD of $M\times N$ matrices. The slower
convergence of WeCAN may be another reason. Fig. \ref{fig:ZCZ_set}
shows the auto- and cross-correlations (normalized by $N$) of the
sequence sets generated by the two algorithms. We can see in Fig.
\ref{fig:ZCZ_set} that the correlation sidelobes of the MM-WeCorr
sequence set are suppressed to almost zero (about -175 dB) at the
required lags, while that of the WeCAN sequence set is much higher.
Another observation is that the cross-correlations at lag $k=0$ for
the WeCAN sequence set are very low, although we did not try to suppress
them. The reason is that in WeCAN, the weight at lag $0$ should be
always positive and in fact large enough to ensure some weight matrix
to be positive semidefinite. Thus the ``0-lag'' correlations are
in fact emphasized the most in WeCAN. Note that in MM-WeCorr, the
weight at lag $0$, i.e., $w_{0}$, can take any nonnegative value,
thus it is more flexible to some extent.

\noindent \begin{center}
\begin{figure*}[tbph]
\noindent \centering{}\includegraphics[width=0.8\paperwidth]{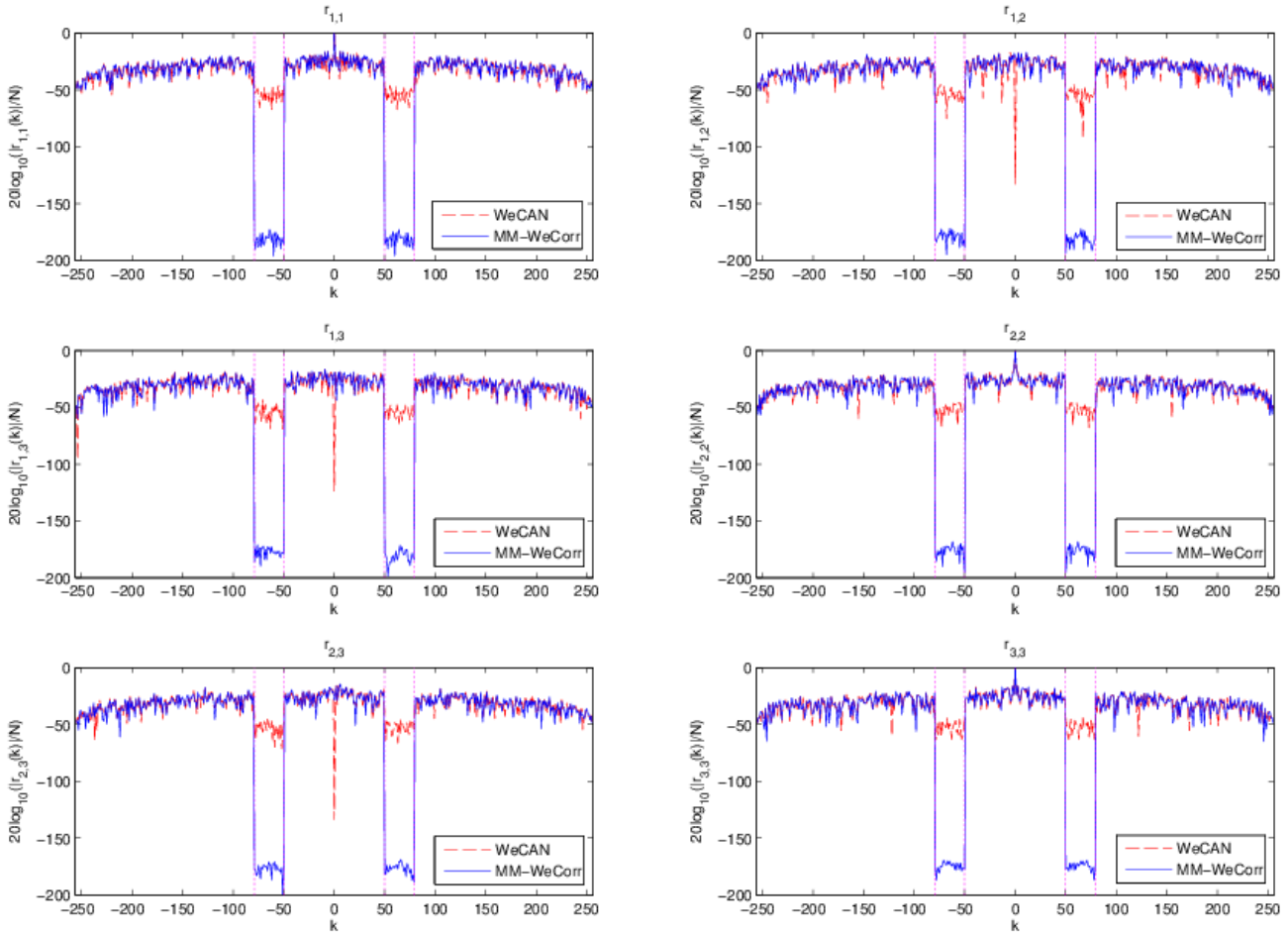}\protect\caption{\label{fig:ZCZ_set}Auto- and cross-correlations of the 256-by-3 sequence
sets generated by MM-WeCorr and WeCAN.}
\end{figure*}

\par\end{center}

\section{Conclusion\label{sec:Conclusion}}

In this paper, we have developed several efficient MM algorithms which
can be used to design unimodular sequence sets with almost complementary
autocorrelations or with both good auto- and cross-correlations. The
proposed algorithms can be viewed as extensions of some single sequence
design algorithms in the literature and share the same convergence
properties, i.e., the convergence to a stationary point. In addition,
all the algorithms can be implemented via FFT and thus are computationally
very efficient. Numerical experiments show that the proposed CSS design
algorithm can generate an almost complementary set of sequences as
long as the cardinality of the set is not too small. In the case of
sequence set design for both good auto- and cross-correlation properties,
the proposed algorithms can get as close to the lower bound of the
correlation criterion as the state-of-the-art method and are much
faster. It has also been observed that the proposed weighted correlation
minimization algorithm can produce sets of unimodular sequences with
virtually zero auto- and cross-correlations at specified time lags.

\appendices{}

\section{Proof of Lemma \ref{lem:eig_upperBound}\label{sec:Proof-of-Lemma-eigUB}}
\begin{IEEEproof}
First, with Lemma \ref{lem:eig_set}, we have 
\begin{equation}
\begin{aligned} & \,\,\lambda_{{\rm max}}\left(\mathbf{R}-\mathbf{B}\circ\big(\mathbf{x}^{(l)}(\mathbf{x}^{(l)})^{H}\big)\right)\\
\leq & \,\,\lambda_{{\rm max}}(\mathbf{R})-\lambda_{{\rm min}}\left(\mathbf{B}\circ\big(\mathbf{x}^{(l)}(\mathbf{x}^{(l)})^{H}\big)\right)\\
= & \,\,\lambda_{{\rm max}}(\mathbf{R})-\lambda_{{\rm min}}\left(\mathbf{B}\right).
\end{aligned}
\end{equation}
Then, according to Lemma \ref{lem:eig_kron}, it is easy to see that
\begin{equation}
\lambda_{{\rm min}}\left(\mathbf{B}\right)=\begin{cases}
\min\{M\lambda_{{\rm min}}(\mathbf{W}),0\}, & M\geq2\\
\lambda_{{\rm min}}(\mathbf{W}), & M=1,
\end{cases}
\end{equation}
and noticing the fact that $\mathbf{W}$ is symmetric Toeplitz, we
know from Lemma \ref{lem:eig_bounds} that 
\begin{equation}
\lambda_{{\rm min}}(\mathbf{W})\geq\lambda_{W}.
\end{equation}
Thus, 
\begin{equation}
\lambda_{{\rm min}}\left(\mathbf{B}\right)\geq\lambda_{B}=\begin{cases}
\min\left\{ M\lambda_{W},0\right\} , & M\geq2\\
\lambda_{W}, & M=1,
\end{cases}
\end{equation}
and we have 
\begin{equation}
\lambda_{{\rm max}}\left(\mathbf{R}-\mathbf{B}\circ\big(\mathbf{x}^{(l)}(\mathbf{x}^{(l)})^{H}\big)\right)\leq\left\Vert \mathbf{R}\right\Vert -\lambda_{B},
\end{equation}
where $\left\Vert \mathbf{R}\right\Vert $ can be any submultiplicative
matrix norm of $\mathbf{R}$.
\end{IEEEproof}

\section{Proof of Lemma \ref{lem:diagonal}\label{sec:Proof-of-Lemma-diag}}
\begin{IEEEproof}
The $N\times N$ Toeplitz matrix $\mathbf{T}$ can be embedded in
a circulant matrix $\mathbf{C}$ of dimension $2N\times2N$ as follows:
\begin{equation}
\mathbf{C}=\left[\begin{array}{cc}
\mathbf{T} & \mathbf{W}\\
\mathbf{W} & \mathbf{T}
\end{array}\right],
\end{equation}
where
\begin{equation}
\mathbf{W}=\begin{bmatrix}0 & t_{1-N} & \cdots & t_{-1}\\
t_{N-1} & 0 & \ddots & \vdots\\
\vdots & \ddots & \ddots & t_{1-N}\\
t_{1} & \cdots & t_{N-1} & 0
\end{bmatrix}.
\end{equation}
The circulant matrix $\mathbf{C}$ can be diagonalized by the FFT
matrix \cite{Gray2006}, i.e.,
\begin{equation}
\mathbf{C}=\frac{1}{2N}\mathbf{F}^{H}{\rm Diag}(\mathbf{F}\mathbf{c})\mathbf{F},\label{eq:circ_decomp}
\end{equation}
where $\mathbf{c}$ is the first column of $\mathbf{C},$ i.e., $\mathbf{c}=[t_{0},t_{-1},\cdots,t_{1-N},0,t_{N-1},\cdots,t_{1}]^{T}.$
Since the matrix $\mathbf{T}$ is just the upper left $N\times N$
block of $\mathbf{C}$, we can easily obtain $\mathbf{T}=\frac{1}{2N}\mathbf{F}_{:,1:N}^{H}{\rm Diag}(\mathbf{F}\mathbf{c})\mathbf{F}_{:,1:N}$.
\end{IEEEproof}

\section{Proof of Lemma \ref{lem:quartic_major} \label{sec:Proof-of-Lemma-quartic-major}}
\begin{IEEEproof}
For any given $x_{0}\in[0,t)$, let us consider the quadratic function
of the following form: 
\begin{equation}
g(x|x_{0})=x_{0}^{4}+4x_{0}^{3}(x-x_{0})+a(x-x_{0})^{2},\label{eq:g_x}
\end{equation}
where $a>0$. It is easy to check that $f(x_{0})=g(x_{0}|x_{0}).$
So to make $g(x|x_{0})$ be a majorization function of $f(x)$ at
$x_{0}$ over the interval $[0,t]$, we need to further have $f(x)\leq g(x|x_{0})$
for all $x\in[0,t]$. Equivalently, we must have
\begin{equation}
\begin{aligned}a & \geq\frac{x^{4}-x_{0}^{4}-4x_{0}^{3}(x-x_{0})}{(x-x_{0})^{2}}\\
 & =x^{2}+2x_{0}x+3x_{0}^{2}
\end{aligned}
\label{eq:a_inequal}
\end{equation}
for all $x\in[0,t]$. Let us define the function 
\begin{equation}
A(x|x_{0})=x^{2}+2x_{0}x+3x_{0}^{2},
\end{equation}
then condition \eqref{eq:a_inequal} is equivalent to 
\begin{equation}
\begin{aligned}a & \geq\max_{x\in[0,t]}\,A(x|x_{0}).\end{aligned}
\label{eq:a_inequal2}
\end{equation}
Since the derivative of $A(x|x_{0})$, given by 
\begin{equation}
A^{\prime}(x|x_{0})=2x+2x_{0},
\end{equation}
is nonnegative for all $x\in[0,t],$ we know that $A(x|x_{0})$ is
nondecreasing on the interval $[0,t]$ and the maximum is achieved
at $x=t.$ Thus, condition \eqref{eq:a_inequal2} becomes 
\begin{equation}
\begin{aligned}a & \geq A(t|x_{0})\\
 & =t^{2}+2x_{0}t+3x_{0}^{2}.
\end{aligned}
\end{equation}
Finally, by appropriately rearranging the terms of $g(x|x_{0})$ in
\eqref{eq:g_x}, we can obtain the function in \eqref{eq:quad_major}.
The proof is complete.
\end{IEEEproof}
\bibliographystyle{IEEEtran}
\bibliography{sequence_set}

\begin{thebibliography}{10}
\providecommand{\url}[1]{#1}
\csname url@samestyle\endcsname
\providecommand{\newblock}{\relax}
\providecommand{\bibinfo}[2]{#2}
\providecommand{\BIBentrySTDinterwordspacing}{\spaceskip=0pt\relax}
\providecommand{\BIBentryALTinterwordstretchfactor}{4}
\providecommand{\BIBentryALTinterwordspacing}{\spaceskip=\fontdimen2\font plus
\BIBentryALTinterwordstretchfactor\fontdimen3\font minus
  \fontdimen4\font\relax}
\providecommand{\BIBforeignlanguage}[2]{{%
\expandafter\ifx\csname l@#1\endcsname\relax
\typeout{** WARNING: IEEEtran.bst: No hyphenation pattern has been}%
\typeout{** loaded for the language `#1'. Using the pattern for}%
\typeout{** the default language instead.}%
\else
\language=\csname l@#1\endcsname
\fi
#2}}
\providecommand{\BIBdecl}{\relax}
\BIBdecl

\bibitem{he2012waveform}
H.~He, J.~Li, and P.~Stoica, \emph{{Waveform Design for Active Sensing Systems:
  A Computational Approach}}.\hskip 1em plus 0.5em minus 0.4em\relax Cambridge
  University Press, 2012.

\bibitem{levanon2004radar}
N.~Levanon and E.~Mozeson, \emph{{Radar Signals}}.\hskip 1em plus 0.5em minus
  0.4em\relax John Wiley \& Sons, 2004.

\bibitem{stoica2009new}
P.~Stoica, H.~He, and J.~Li, ``New algorithms for designing unimodular
  sequences with good correlation properties,'' \emph{IEEE Transactions on
  Signal Processing}, vol.~57, no.~4, pp. 1415--1425, 2009.

\bibitem{MISL}
J.~Song, P.~Babu, and D.~P. Palomar, ``Optimization methods for designing
  sequences with low autocorrelation sidelobes,'' \emph{IEEE Transactions on
  Signal Processing}, vol.~63, no.~15, pp. 3998--4009, Aug. 2015.

\bibitem{WISL_song2015}
\BIBentryALTinterwordspacing
------, ``Sequence design to minimize the weighted integrated and peak sidelobe
  levels,'' \emph{submitted to IEEE Transactions on Signal Processing}, 2015.
  [Online]. Available: \url{http://arxiv.org/abs/1506.04234}
\BIBentrySTDinterwordspacing

\bibitem{compleSet_Searle2008}
S.~Searle, S.~Howard, and B.~Moran, ``The use of complementary sets in {MIMO}
  radar,'' in \emph{2008 42nd Asilomar Conference on Signals, Systems and
  Computers}, Oct. 2008, pp. 510--514.

\bibitem{compleSet_Levanon2009}
N.~Levanon, ``Noncoherent radar pulse compression based on complementary
  sequences,'' \emph{IEEE Transactions on Aerospace and Electronic Systems},
  vol.~45, no.~2, pp. 742--747, April 2009.

\bibitem{compleSet_Schmidt2007}
K.~Schmidt, ``Complementary sets, generalized {Reed-Muller} codes, and power
  control for {OFDM},'' \emph{IEEE Transactions on Information Theory},
  vol.~53, no.~2, pp. 808--814, Feb. 2007.

\bibitem{compleSet_Garcia2010}
E.~Garcia, J.~Garcia, J.~Urena, M.~Perez, and D.~Ruiz, ``Multilevel
  complementary sets of sequences and their application in {UWB},'' in
  \emph{2010 International Conference on Indoor Positioning and Indoor
  Navigation (IPIN)}, Sep. 2010, pp. 1--5.

\bibitem{compleSet_Tseng2000}
S.-M. Tseng and M.~Bell, ``Asynchronous multicarrier {DS-CDMA} using mutually
  orthogonal complementary sets of sequences,'' \emph{IEEE Transactions on
  Communications}, vol.~48, no.~1, pp. 53--59, Jan. 2000.

\bibitem{compleSet_Spasojevic2001}
P.~Spasojevic and C.~Georghiades, ``Complementary sequences for {ISI} channel
  estimation,'' \emph{IEEE Transactions on Information Theory}, vol.~47, no.~3,
  pp. 1145--1152, Mar. 2001.

\bibitem{SequenceSet_Mojtaba2013}
M.~Soltanalian, M.~M. Naghsh, and P.~Stoica, ``A fast algorithm for designing
  complementary sets of sequences,'' \emph{Signal Processing}, vol.~93, no.~7,
  pp. 2096--2102, 2013.

\bibitem{lowCorrSet_Oppermann1997}
I.~Oppermann and B.~Vucetic, ``Complex spreading sequences with a wide range of
  correlation properties,'' \emph{IEEE Transactions on Communications},
  vol.~45, no.~3, pp. 365--375, Mar. 1997.

\bibitem{SequenceSet_He2009}
H.~He, P.~Stoica, and J.~Li, ``Designing unimodular sequence sets with good
  correlations--including an application to {MIMO} radar,'' \emph{IEEE
  Transactions on Signal Processing}, vol.~57, no.~11, pp. 4391--4405, Nov.
  2009.

\bibitem{hunter2004MMtutorial}
D.~R. Hunter and K.~Lange, ``A tutorial on {MM} algorithms,'' \emph{The
  American Statistician}, vol.~58, no.~1, pp. 30--37, 2004.

\bibitem{MM_Stoica}
P.~Stoica and Y.~Selen, ``Cyclic minimizers, majorization techniques, and the
  expectation-maximization algorithm: a refresher,'' \emph{IEEE Signal
  Processing Magazine}, vol.~21, no.~1, pp. 112--114, Jan. 2004.

\bibitem{razaviyayn2013unified}
M.~Razaviyayn, M.~Hong, and Z.-Q. Luo, ``A unified convergence analysis of
  block successive minimization methods for nonsmooth optimization,''
  \emph{SIAM Journal on Optimization}, vol.~23, no.~2, pp. 1126--1153, 2013.

\bibitem{roger1994topics}
H.~Roger and R.~J. Charles, \emph{{Topics in Matrix Analysis}}.\hskip 1em plus
  0.5em minus 0.4em\relax Cambridge University Press, 1994.

\bibitem{eig_localization}
P.~J. S.~G. Ferreira, ``Localization of the eigenvalues of {Toeplitz} matrices
  using additive decomposition, embedding in circulants, and the {Fourier}
  transform,'' in \emph{Proceedings of the 10th IFAC Symposium on System
  Identification}, Copenhagen, Denmark, Jul. 1994, pp. 271--276.

\bibitem{Bertsekas2003}
D.~P. Bertsekas, A.~Nedi\'{c}, and A.~E. Ozdaglar, \emph{{Convex Analysis and
  Optimization}}.\hskip 1em plus 0.5em minus 0.4em\relax Athena Scientific,
  2003.

\bibitem{SQUAREM}
R.~Varadhan and C.~Roland, ``Simple and globally convergent methods for
  accelerating the convergence of any {EM} algorithm,'' \emph{Scandinavian
  Journal of Statistics}, vol.~35, no.~2, pp. 335--353, 2008.

\bibitem{Gray2006}
R.~M. Gray, \emph{{Toeplitz and Circulant Matrices: A review}}.\hskip 1em plus
  0.5em minus 0.4em\relax Now Publishers Inc, Jan. 2006, vol.~2, no.~3.

\end{thebibliography}

\end{document}